\documentclass[12pt,english]{article}
\usepackage[cp1251]{inputenc}
\usepackage[english]{babel}
\usepackage[T1]{fontenc}
\usepackage{amsthm, amsfonts,amsmath,amssymb,graphicx, mathrsfs,url}
\usepackage{rotating}

\begin{document}

\newtheorem{teo}{Theorem}
\newtheorem*{teon}{Theorem}
\newtheorem{lem}{Lemma}
\newtheorem*{lemn}{Lemma}
\newtheorem{prp}{Proposition}
\newtheorem*{prpn}{Proposition}
\newtheorem{ass}{Assertion}
\newtheorem*{assn}{Assertion}
\newtheorem{assum}{Assumption}
\newtheorem*{assumn}{Assumption}
\newtheorem{stat}{Statement}
\newtheorem*{statn}{Statement}
\newtheorem{cor}{Corollary}
\newtheorem*{corn}{Corollary}
\newtheorem{hyp}{Hypothesis}
\newtheorem*{hypn}{Hypothesis}
\newtheorem{con}{Conjecture}
\newtheorem*{conn}{Conjecture}
\newtheorem{dfn}{Definition}
\newtheorem*{dfnn}{Definition}
\newtheorem{problem}{Problem}
\newtheorem*{problemn}{Problem}
\newtheorem{notat}{Notation}
\newtheorem*{notatn}{Notation}
\newtheorem{quest}{Question}
\newtheorem*{questn}{Question}

\newtheorem{rem}{Remark}
\newtheorem*{remn}{Remark}
\newtheorem{exa}{Example}
\newtheorem*{exan}{Example}
\newtheorem{cas}{Case}
\newtheorem*{casn}{Case}
\newtheorem{claim}{Claim}
\newtheorem*{claimn}{Claim}
\newtheorem{com}{Comment}
\newtheorem*{comn}{Comment}


\newtheorem*{proofn}{Proof}

\selectlanguage{english}
%

\def\JournalNumber{0}
\def\JournalVolume{00}
%



\def\Z{\mathbb Z}
\def\R{\mathbb R}
\newcommand{\A}{\operatorname{\mathcal{A}}\nolimits}
\newcommand{\Id}{\operatorname{Id}\nolimits}
\newcommand{\spann}{\operatorname{span}\nolimits}
\newcommand{\Conj}{\operatorname{Conj}\nolimits}
\newcommand{\Max}{\operatorname{Max}\nolimits}
\newcommand{\cl}{\operatorname{cl}\nolimits}
\newcommand{\restr}[2]{\left. #1 \right|_{#2}}


\title{Maxwell Strata and Cut Locus in Sub-Riemannian Problem on Engel group}
\author{Andrey\,A.~Ardentov, Yuri\,L.~Sachkov \\ \small{aaa@pereslavl.ru, yusachkov@gmail.com}  \\ \small{Program Systems Institute of RAS} \\ \footnotesize{Pereslavl-Zalessky, Yaroslavl Region, 152020, Russia}}
\maketitle

\begin{abstract}
We consider the nilpotent left-invariant sub-Riemannian structure on the Engel group. This structure gives a fundamental local approximation of a generic rank 2 sub-Riemannian structure on a 4-manifold near a generic point (in particular, of the kinematic models of a car with a trailer). On the other hand, this is the simplest
 sub-Riemannian structure of step three.

We describe the global structure of the cut locus (the set of points where geodesics lose their global optimality), the Maxwell set (the set of points that admit more than one minimizer), and the intersection of the cut locus with the caustic (the set of conjugate points along all geodesics). 

The group of symmetries of the cut locus is described: it is generated by
a one-parameter group of dilations $\R_+$ and a discrete group of reflections $\Z_2 \times \Z_2 \times \Z_2$.

The cut locus admits a stratification with 6 three-dimensional strata, 12 two-dimensional strata, and 2 one-dimensional strata.
Three-dimen-sional strata of the cut locus are Maxwell strata of multiplicity 2 (for each point there are 2 minimizers). Two-dimensional strata of the cut locus consist of conjugate points. Finally, one-dimensional strata are Maxwell strata of infinite multiplicity, they consist of conjugate points as well.

Projections of sub-Riemannian geodesics to the 2-dimensional plane of the distribution are Euler elasticae.
For each point of the cut locus, we describe the Euler elasticae corresponding to minimizers coming to this point.

Finally, we describe the structure of the optimal synthesis, i.e., the set of minimizers for each terminal point in the Engel group.

\end{abstract}



\def\ds{\displaystyle}
\def\cut{\operatorname{cut}}
\def\tcut{t_{\cut}}
\def\Cut{\operatorname{Cut}}
\newcommand{\vect}[1]{\left( \begin{array}{c} #1 \end{array} \right)}
\newcommand{\Exp}{\operatorname{Exp}\nolimits}
\newcommand{\sgn}{\operatorname{sgn}\nolimits}
\newcommand{\E}{\operatorname{E}}
\newcommand{\sn}{\operatorname{sn}\nolimits}
\newcommand{\cn}{\operatorname{cn}\nolimits}
\newcommand{\dn}{\operatorname{dn}\nolimits}
\newcommand{\am}{\operatorname{am}\nolimits}
\newcommand{\FIX}{\operatorname{FIX}\nolimits}
\newcommand{\MAX}{\operatorname{MAX}\nolimits}
\def\newtilde{{\,\,\approx\,\,}}
\newcommand{\CMAX}{\operatorname{CMAX}\nolimits}
\def\conj{\operatorname{conj}}
\newcommand{\I}{\operatorname{\mathcal{I}}\nolimits}
\newcommand{\N}{\operatorname{\mathcal{N}}\nolimits}
\newcommand{\Fix}{\operatorname{Fix}\nolimits}
\newcommand{\fix}{\operatorname{fix}\nolimits}


\begin{flushright}
\textit{To the memory of Vladimir Arnold}
\end{flushright}

\section{INTRODUCTION}
The Engel group $M$ is a nilpotent four-dimensional Lie group, connected and simply connected, which has Lie algebra 
 $L =  \spann(X_1, X_2, X_3, X_4)$ with the multiplication table
$$
[X_1, X_2] = X_3, \qquad
[X_1, X_3] = X_4, \qquad 
[X_2,X_3]=[X_1,X_4]=[X_2,X_4]=0. 
$$

This article studies the sub-Riemannian structure~\cite{GCT} on the Lie group $M$  generated by the left-invariant orthonormal frame $X_1$, $X_2$. This structure gives a nilpotent approximation~\cite{ABB} to a generic rank two sub-Riemannian structure in four-dimensional space near a generic point. 

In certain coordinates $(x,y,z,w)$ on the Engel group $M \cong \R^4$, the nilpotent sub-Riemannian problem is stated as follows:
\begin{align*}
&\left\{
\begin{array} {l}
\dot x = u_1, \qquad
\dot y = u_2, \qquad
\dot z = -u_1 \frac{y}{2} + u_2 \frac{x}{2}, \qquad
\dot w = u_2 \frac{x^2}{2}, \\
q = (x,y,z,w) \in {\mathbb R}^4, \quad (u_1, u_2) \in  {\mathbb R}^2,\\
q(0) = q_0 = (0,0,0,0), \quad q(t_1) = (x_1, y_1, z_1, w_1),  
\end{array}\right. \\
&l = \int_0^{t_1} \sqrt{u_1^2 + u_2^2} \, d t \rightarrow \min. 
\end{align*}

This paper has the following structure. 
In Subsection~\ref{subsec:previous} we recall results on the problem obtained in previous works~\cite{engel, engel_conj, engel_cut}.
In Subsection~\ref{subsec:preliminary} we prove some simple preliminary results on the cut locus.
In Sections~\ref{sec:x=z=0}, \ref{sec:z=0+}, \ref{sec:z=0}, \ref{sec:x=0+},  \ref{sec:x=0} we describe respectively intersection of the cut locus with the sets $\{x=z=0\}$, $\{z=0, \ x > 0\}$,  $\{z=0\}$, $\{x=0, \ z > 0\}$, and $\{x=0\}$. In Section~\ref{sec:cut_glob} we sum up these results by describing a global stratification of the cut locus. In Appendices A--C we prove technical lemmas from Sections~\ref{sec:x=z=0}, \ref{sec:z=0+}, and \ref{sec:x=0+}.

\subsection{Previously obtained results}\label{subsec:previous}

This paper continues a series of works~\cite{engel, engel_conj, engel_cut}, where a detailed study of the sub-Riemannian problem on the Engel group was started (in these works, instead of the coordinate $w$, we used the coordinate $v = w + y^3/6$). First we recall the main results of these works. 

For each point $q_1 \in M$ there exists an optimal trajectory (sub-Riemannian minimizer). Sub-Riemannian geodesics are described by Pontryagin maximum principle. 
Abnormal trajectories are simultaneously normal, and their endpoints fill two rays
\begin{equation*}
\A_{\pm} = \{q \in M \mid x = z = w = 0, \sgn y = \pm 1\}.
\end{equation*}

Geodesics are parametrized by Jacobi elliptic functions~\cite{whit_watson}. Projections of geodesics onto the plane  $(x,y)$ are Euler's elasticae~\cite{euler, love}. Small arcs of geodesics are optimal; however, large arcs are, in general, not optimal. A point at which a geodesic loses its optimality is called a cut point. The union of cut points along all geodesics is called the cut locus. The cut locus is one of the most important characteristics of sub-Riemannian structures~\cite{ABB}, and the main goal of this paper is its description for the left-invariant sub-Riemannian problem on the Engel group.

A generic geodesic loses its optimality at a Maxwell point, i.e., a point where several geodesics of the same length meet one another. Maxwell points are fixed points of discrete symmetries $\varepsilon^1$, $\varepsilon^2$, $\varepsilon^4$  acting on the Engel group as follows:
\begin{align}
\varepsilon^1(x,y,z,w) &= (x,y,-z,w-xz), \label{vareps1}\\
\varepsilon^2(x,y,z,w) &= (-x,y,z,w-xz), \label{vareps2}\\
\varepsilon^4(x,y,z,w) &= (-x,-y,z,-w). \label{vareps4}
\end{align}
The subspaces  
$
M_x = \{ q \in M \mid x = 0\}$, 
$M_z = \{ q \in M \mid z = 0\}$
 consist of fixed points of symmetries~(\ref{vareps1}), (\ref{vareps2}), the cut locus is contained in the union of these subspaces (see further Th.~\ref{th:cut_simple}).

The problem has also continuous symmetries --- the one-parameter group of dilations given by the flow of the vector field
$$  X_0 = x \frac{\partial }{\partial x} + y \frac{\partial }{\partial y} + 2 z \frac{\partial }{\partial z} + 3w \frac{\partial }{\partial w} .$$

The main result of works~\cite{engel, engel_conj, engel_cut} is the explicit description of the cut time (see Th.~\ref{tcut} below) and the proof that there is a unique minimizer for each point $q_1 \in M \cap \{ x z \neq 0\}$ (see Th.~\ref{NM} below). 

We showed in work~\cite{engel} that the family of all extremal trajectories of the problem is parametrized by the cylinder  
\begin{align*}
C = \left\{\lambda \in T_{q_0} ^* M \mid H(\lambda)= 1/2\right\}= \left\{\lambda = (\theta, c, \alpha  ) \mid \theta \in S^1, \ c, \alpha \in \R  \right\},
\end{align*}
where $H = (\langle \lambda, X_1\rangle^2 + \langle \lambda, X_2\rangle^2)/2$, $\lambda \in T^*M$,  is the maximized Hamiltonian of Pontryagin maximum principle; $(\theta, c, \alpha)$ are certain natural coordinates on the cylinder  $C$. Parametrization of extremal trajectories is defined by the exponential mapping 
\begin{align*}
&\Exp \colon N \to M,  &&N = C \times \R_+,\\
&\Exp (\nu) = q_t = (x_t, y_t, z_t, w_t),  &&\nu = (\lambda, t).
\end{align*}
The function $E =  c^2/2 - \alpha \cos \theta$ is constant along extremal trajectories, and the cylinder $C$ stratifies according to its values: 
\begin{align*}
C&=\sqcup_{i=1}^7 C_i,    \  &&C_4 = \big\{\lambda \in C \mid \alpha \neq 0, E=-|\alpha|\big\}, 
\\   
C_1 &= \big\{\lambda \in C \mid \alpha \neq 0, E\in(- |\alpha|, 
|\alpha|)\big\}, \  &&C_5 = \big\{\lambda \in C \mid \alpha \neq 0, E=|\alpha|, c = 0\big\}, 
\\
C_2 &= \big\{\lambda \in C \mid \alpha \neq 0, E\in(|\alpha|,+\infty)\big\}, \  &&C_6 = \big\{\lambda \in C \mid \alpha = 0, \ c \neq 0\big\},  
\\
C_3 &= \big\{\lambda \in C \mid \alpha \neq 0, E=|\alpha|, c \neq 0 \big\}, \  &&C_7 = \big\{\lambda \in C \mid \alpha = c = 0\big\}. 
\end{align*}
Denote the corresponding subsets in preimage of the exponential mapping: 
$N_i = C_i \times \R_+$.

An arbitrary extremal trajectory of the sub-Riemannian problem on the Engel group projects to the plane 
 $(x,y)$ into an Euler elastica. Each subset   $C_i, \ i=1, \dots, 7,$ corresponds to a certain type of Euler elasticae. 

For parametrization of trajectories for the subsets  
$C_1$ (inflectional elasticae), 
$C_2$ (noninflectional elasticae) and
$C_3$ (critical elasticae) were introduced elliptic coordinates
$\lambda = (\varphi, k, \alpha)$. The parameter $k$ is a reparametrization of the first integral $E$. On the set $C_3$ we have $k=1$, this set separates  $C_1$ and $C_2$, where $k\in(0,1)$.  The remaining subsets
$C_i$, $i = 4, \dots, 7,$ (lines and circles in the plane   $(x,y)$) are parametrized by the coordinates  $\lambda = (\theta,c, \alpha)$. Notice that trajectories for the case   $\lambda \in C_4 \cup C_5$ are defined by formulas for   $\lambda \in C_7$, when $\sin \theta = 0$, thus the cases   $C_4$, $C_5$ are not considered in this paper.

On the subsets $N_1, N_2, N_3$  were introduced new coordinates $(p,\tau,\sigma)$:
\begin{align*}
\sigma &= \sgn \alpha \sqrt{ |\alpha| },&&&& \\
(\lambda,t) &\in N_1 \cup N_3 &\Rightarrow \quad\qquad &p = |\sigma| t/2, && \tau = |\sigma| (\varphi + t/2),\\
(\lambda,t) &\in N_2 &\Rightarrow \quad\qquad  &p = |\sigma| t/(2 k), &&\tau = |\sigma| (\varphi + t/2)/k.
\end{align*}

\begin{dfn}
The cut time $\tcut(\lambda)$ is the time when the extremal trajectory corresponding to the covector   $\lambda$ loses its global optimality:
$$
\tcut(\lambda) = \sup \big\{ t>0 \mid \Exp (\lambda,s) \text{ optimal for } s\in[0,t]\big\}.
$$

The cut locus is the set
$
\Cut = \big\{\Exp(\lambda,t) \mid \lambda \in C, t = \tcut (\lambda)\big\}.
$

\end{dfn}
\begin{rem}
If $\tcut(\lambda) = \infty$, then the trajectory $\Exp (\lambda,s)$ is optimal on the whole ray    $s\in[0,\infty)$. In this case the nonstrict inequality   $s\leq\tcut(\lambda)$ should be understood as a strict one.
\end{rem}

\begin{teo}[\cite{engel_cut}~Cor.~4.2]\label{tcut} The cut time has the following explicit expression: 
\begin{align*}
& \forall \lambda \in C_1 &\quad & \tcut (\lambda) = \frac{\min \big(2 p_z^1(k), 4 K(k)\big)}{|\sigma|} = \left\{\begin{array}{l}
			4 K(k)/|\sigma|,  \  k\in (0,k_0], \\
			2 p_z^1(k)/|\sigma|, \ k\in[k_0,1),
			\end{array}\right. \\
& \forall \lambda \in C_2 &\quad & \tcut (\lambda) = \frac{2 K(k) k}{|\sigma|}, \\
& \forall \lambda \in C_6 &\quad & \tcut (\lambda) = \frac{2 \pi}{ \sqrt{|c| } }, \\
& \forall \lambda \in C_3 \cup C_7 &\quad & \tcut (\lambda) = +\infty.
\end{align*}
\end{teo}
Here $K(k) = \int_0^\frac{\pi}{2} \frac{dt}{\sqrt{1- k^2 \sin^2 t}}$ is the complete elliptic integral of the first kind; $p_z^1(k)\in \big(K(k), 3 K(k)\big)$ is the first positive root of the function   $f_z(p,k)=\dn p \,\sn p+ \big(p-2\E(p)\big)\cn p.$ The functions $\sn p, \cn p, \dn p$ are Jacobian elliptic functions with modulus   $k$ by default (since the modulus  $k$ is constant along extremal trajectories), i.e., e.g.   $\sn p = \sn (p,k)$; also $\ds \E(p)=\int_0^p \dn^2 t \,dt.$

On the subsets $N_1$, $N_2$ the coordinates $(p, \tau)$ are transformed respectively into   $(u_1, u_2)$ via the formulas
$
u_1 = \am p$,   $u_2= \am \tau$,
where $\am$ is the elliptic amplitude, inverse function to the incomplete elliptic integral of the first kind:  $F (\am p) = p$. Summing up, we use the following coordinates for parametrization of the exponential mapping on subsets: 
\begin{align}
&\nu \in N_1 \cup N_2, && \nu = (k, u_1, u_2, \sigma), \label{cord1} \\
&\nu \in N_3, && \nu = (p,\tau,\sigma), \\
&\nu \in N_6, && \nu = (\theta, c, t), \\
&\nu \in N_7, && \nu = (\theta, t). \label{cord4}
\end{align}
Explicit formulas for the exponential mapping 
  $\Exp(\lambda, t)$ for each subset   $N_i$, $i=1,\dots,7,$ are given in paper~\cite{engel_cut}.

In work~\cite{engel} we showed that vanishing one of the coordinates   $x$, $z$ on geodesics is related to discrete symmetries   $\varepsilon^1$--$\varepsilon^7$  of the exponential mapping. The group of discrete symmetries   $\{\Id, \varepsilon^1, \dots, \varepsilon^7\}$ is isomorphic to the group of symmetries of parallelepiped   $\Z_2 \times \Z_2 \times\Z_2$. Explicit expression of the symmetries 
  $\varepsilon^1$, $\varepsilon^2$, $\varepsilon^4$ is presented above in~(\ref{vareps1})--(\ref{vareps4}).
	The sets
  $M_z$, $M_x$ consist of fixed points of symmetries~(\ref{vareps1}), (\ref{vareps2}) in the group  $M$. Action of the symmetries    $\varepsilon^1$, $\varepsilon^2$ in the preimage $N$ is expressed as follows:
\begin{align}
\varepsilon^1(\theta, c, \alpha, t) &= (\theta_t, -c_t, \alpha, t), \label{vareps1a}\\
\varepsilon^2(\theta, c, \alpha, t) &= (-\theta_t, c_t, \alpha, t), \label{vareps2a}\\
\varepsilon^4(\theta, c, \alpha, t) &= (\theta + \pi, c, -\alpha, t). \label{vareps4a}
\end{align}

Conditions of invariance of points in the preimage of the exponential mapping w.r.t. the symmetries
  $\varepsilon^1$, $\varepsilon^2$, are given respectively by the equalities   $c_{t/2} = 0$ and $\sin \theta_{t/2} = 0$. On the subsets   $N_i, i=1,\dots,7,$ in coordinates~(\ref{cord1})--(\ref{cord4}) these conditions are described by Table~\ref{tab:fix}.

\begin{table}[!ht]
\caption{\label{tab:fix} Conditions of invariance of points w.r.t. action of   $\varepsilon^1, \varepsilon^2$ in $N$}
\begin{tabular}{|c|c|c|c|c|c|}
\hline
$\nu \in N_i$ & $N_1$ & $N_2$ & $N_3$ & $N_6$ & $N_7$ \\
\hline
$\varepsilon^1(\nu) = \nu$ &  $\cos u_2 = 0$ & $\emptyset$ & $\emptyset$ & $\emptyset$ & $N_7$ \\
\hline
$\varepsilon^2(\nu) = \nu$ &  $\sin u_2 = 0$ & $\sin u_2 \cos u_2 = 0$ & $\tau = 0$ & $2\theta + c t = \underset{n\in \Z}{2 \pi n,}$ & $\sin \theta =0$ \\
\hline
\end{tabular}
\end{table}

In the preimage of the exponential mapping, optimal trajectories correspond to the set
\begin{align*}
\widehat{N} &= \big\{(\lambda,t) \in N \mid t \leq \tcut (\lambda) \big\}.
\end{align*}

After having excluded the initial point  
  $q_0$, we get the following set of terminal points  $q_1$:
\begin{align*}
\widehat{M} &= \big\{(x,y,z,w) \in \R^4 \mid x^2+y^2+z^2+w^2 \neq 0 \big\}. 
\end{align*}
Below we mean that 
  $q_1 = (x,y,z,w) \in \widehat{M}$.

According to the action of the main symmetries 
  $\varepsilon^1, \varepsilon^2$, the subsets   $\widehat{M}, \widehat{N}$ decompose into the following subsets:
\begin{align*}
\widehat{M} &= M' \sqcup \widetilde{M},  \qquad \widehat{N} = N' \sqcup \widetilde{N}, \\
M' &= \big\{q \in M \mid x z = 0, x^2 + y^2 + z^2 + w^2 \neq 0\big\}, \quad \widetilde{M} = \big\{q \in M \mid x z \neq 0\big\}, \\
N' &= \big\{(\lambda,t) \in N \mid t = \tcut (\lambda) \text{ or } c_{t/2} \sin \theta_{t/2} = 0, t < \tcut (\lambda)\big\}, \\
\widetilde{N} &= \big\{(\lambda,t) \in N \mid t < \tcut (\lambda), c_{t/2} \sin \theta_{t/2} \neq 0\big\}. 
\end{align*}

\begin{teo}[\cite{engel_cut}~Cor.~3.21] \label{NM}
The mapping $\Exp\colon \widetilde{N} \to \widetilde{M}$ is a diffeomorphism.
\end{teo}

\subsection{Preliminary results on cut locus}\label{subsec:preliminary}
The set $N'$ consists of cut points and fixed points of the symmetries   $\varepsilon^1, \varepsilon^2$:
\begin{align*}
N' &= N_{\cut} \sqcup \FIX, \qquad \FIX = \FIX^1 \sqcup \FIX^2 \sqcup \FIX^{12}, \\
N_{\cut}  &= \big\{(\lambda,t) \in N \mid t = \tcut(\lambda)\big\}, \qquad \FIX^{12}  = \big\{\nu \in N_7 \mid \sin \theta = 0\big\}, \\
\FIX^i &= \big\{(\lambda, t) \in N \backslash \FIX^{12} \mid t<\tcut(\lambda), \ \varepsilon^i(\lambda,t) = (\lambda, t)\big\}, \quad i=1,2.
\end{align*}
Obviously, $\Cut = \Exp(N_{\cut})$.

The set $M'$ stratifies as follows:
\begin{align*}
M'&= M_{00} \sqcup M_{0+} \sqcup M_{0-} \sqcup M_{+0} \sqcup M_{-0}, \  &M_{00} &= \{q \in M \mid x = z = 0, y^2+w^2 \neq 0 \}, \\
M_{0\pm} &= \{q \in M \mid x=0, \sgn z =\pm 1\}, \  &M_{\pm0} &= \{q \in M \mid \sgn x =\pm 1, z=0\}.
\end{align*}

Certain simple geometric properties of the cut locus follow immediately from previously obtained results.

\begin{teo}\label{th:cut_simple}
\begin{itemize}
\item[$1)$]
$\Cut \subset M_x \cup M_z$,
\item[$2)$]
$\varepsilon^i(\Cut) = \Cut$, $i = 1, \dots, 7$,
\item[$3)$]
$e^{t X_0}(\Cut) = \Cut$, $t\in\R$.
\end{itemize}
\end{teo}
\begin{proofn}
Inclusion $1)$ follows from Th.~\ref{tcut}.
Equalities $2)$, $3)$ follow from invariance of the cut time w.r.t. reflections $\varepsilon^i$ and its homogeneity w.r.t. dilations $e^{tX_0}$, see item (3) of Cor.~4.2~\cite{engel_cut}. 
\end{proofn}

Denote $N_- = \big\{ (\lambda, t) \in N \mid t < \tcut (\lambda) \big\}$.

\begin{lem} \label{cutfix}
We have 
$\Cut \cap \Exp (N_-) = \emptyset,$ thus $\Cut \cap \Exp (\FIX) = \emptyset$.
\end{lem}
\begin{proofn}
By contradiction, let $\Cut \cap \Exp(N_-) \ni q_1 = \Exp(\lambda_1,t_1) = \Exp(\lambda_2, t_2)$,
where $(\lambda_1, t_1) \in N_{\cut}$, $(\lambda_2, t_2) \in N_-$. One can easily see that   $\lambda_1 \neq \lambda_2$ and $q_1(t) = \Exp (\lambda_1, t) \not\equiv \Exp(\lambda_2, t) = q_2(t)$. The both trajectories $q_1(t), t\in [0,t_1],$ and $q_2(t), t \in [0,t_2],$ are optimal, thus    $t_1 = t_2$, and $q_1  = q_2(t_2)$ is a Maxwell point. Thus the trajectory   $q_2 (t)$ is not optimal for   $t>t_2$, which contradicts the inequality   $t_2 <\tcut(\lambda)$.

The second inequality in the statement of this lemma follows from the inclusion
  $\FIX \subset N_-$.
\end{proofn}

\begin{lem}\label{lem0}
Let $\nu \in N$ and $\Exp(\nu) = (x,y,z,w)$. If $\nu \in N_2 \cup N_3 \cup N_6$, then  $z \neq 0$. 
If $\nu \in N_7$, then $z=0$.
\end{lem}
\begin{proofn}
The statement follows immediately from the formulas of the exponential mapping (see~\cite{engel}, Subsecs.~5.5--5.6, 7.4) for the sets $C_i, \ i=2,\dots,7$.
\end{proofn}

\section{Structure of $\Cut \cap \{x=z=0\}$} \label{sec:x=z=0}
Now we study the set
  $N_{00} = \Exp^{-1}(M_{00}) \cap \widehat{N}$. Recall~\cite{engel} that components of geodesics   $x$, $z$ vanish either at the cut time (see Th.~\ref{tcut}) or at a fixed point of a symmetry in the image of the exponential mapping (see Table~\ref{tab:fix}). Denote the function   $u_{1z}(k) = \am \big(p_z^1 (k)\big)$, recall also that   $k_0 \approx 0.909$ is a unique solution to the equation   $2 E(k) - K(k) = 0$, where $E(k)$ is the complete elliptic integral of the second kind.  
	
	We show below in Lemma~\ref{deco0} that points from  $N_{00}$ related to the cut locus belong to the subset   $N_1$ parametrized by the coordinates   $(k, u_1, u_2, \sigma)$. We use these coordinates to define the components   $\MAX_{\pm}^{12}, \MAX_{ij}^{20},\MAX_{ij}^{10}, \ i,j\in\{+,-\}$ from $N_{00}$ related to the cut locus, see Table~\ref{tab:max0}. This table should be read by columns. For example, the first column means that 
\begin{align*}
\MAX_{++}^{20} = \big\{(k, u_1, u_2, \sigma) \in N_1 \mid k \in (0, k_0), u_1 = \pi, u_2 = \frac{\pi}{2}, \sigma \in (0, + \infty) \big\}.
\end{align*}

\begin{table}[ht]
\centering
\caption{\label{tab:max0} Components of $N_{00} \cap N_{\cut} \cap N_1$}
\begin{tabular}{|c|c|c|c|c|c|c|c|c|c|c|c|}
\hline
&\multicolumn{4}{c|}{$\MAX^{20}$}& \multicolumn{2}{c|}{$\MAX^{12}$}& \multicolumn{4}{c|}{$\MAX^{10}$} \\ \cline{2-11}
\begin{sideways}$\MAX$\end{sideways} & \begin{sideways}$\MAX_{++}^{20}$\end{sideways} & \begin{sideways}$\MAX_{+-}^{20}$\end{sideways} & \begin{sideways}$\MAX_{-+}^{20}$\end{sideways} & \begin{sideways}$\MAX_{--}^{20}$\end{sideways} & \begin{sideways}$\MAX_{-}^{12}$\end{sideways} & \begin{sideways}$\MAX_{+}^{12}$\end{sideways} & \begin{sideways}$\MAX_{++}^{10}$\end{sideways} & \begin{sideways}$\MAX_{+-}^{10}$\end{sideways} & \begin{sideways}$\MAX_{-+}^{10}$\end{sideways} & \begin{sideways}$\MAX_{--}^{10}$\end{sideways} \\
\hline
$k$ &  \multicolumn{4}{c|}{$(0,k_0)$} & \multicolumn{2}{c|}{$k_0$} & \multicolumn{4}{c|}{$(k_0,1)$}  \\
\hline
$u_1$ & \multicolumn{4}{c|}{$\pi$} & \multicolumn{2}{c|}{$u_{1z}(k_0)=\pi$} & \multicolumn{4}{c|}{$u_{1z}(k)$} \\
\hline
$u_2$ & $\pi/2$ & $3 \pi/2$ & $\ds \pi/2$ & $3 \pi/2$ & \multicolumn{2}{c|}{$[0, 2 \pi)$} & $0$ & $\pi$ & $0$ & $\pi$ \\
\hline
$\sigma$ & \multicolumn{2}{c|}{$(0,+\infty)$} & \multicolumn{3}{c|}{\,~~~\quad$(-\infty,0)$ \quad~~~}  & \multicolumn{3}{c|}{$(0,+\infty)$} & \multicolumn{2}{c|}{$(-\infty,0)$}  \\
\hline
\end{tabular}
\end{table} 
Let
\begin{align*}
\MAX^{20} &= \MAX_{++}^{20} \sqcup \MAX_{+-}^{20}  \sqcup \MAX_{-+}^{20}  \sqcup \MAX_{--}^{20}, \\
\MAX^{12} &= \MAX_{+}^{12} \sqcup \MAX_{-}^{12}, \\
\MAX^{10} &= \MAX_{++}^{10} \sqcup \MAX_{+-}^{10}  \sqcup \MAX_{-+}^{10}  \sqcup \MAX_{--}^{10}.
\end{align*}

\begin{lem}\label{deco0}
There holds the equality
\begin{align}
N_{00} = \MAX^{20} \sqcup \MAX^{12} \sqcup \MAX^{10} \sqcup \FIX^{12}. \label{N00m}
\end{align}
\end{lem}
\begin{proofn}
Let $\nu = (\lambda,t) \in N_{00}$ and $t<\tcut (\lambda)$. Then $(\lambda,t)$ is a fixed point for the both symmetries   $\varepsilon^1, \varepsilon^2$. One can see from Table~\ref{tab:fix} that the invariance condition holds for the both symmetries only in the case   $\nu \in N_7$,   $\sin \theta = 0$. Further, since we have   $z=0$ in the image, then Lemma~\ref{lem0} implies the inclusion  
$N_{00} \subseteq  \{(\lambda,t) \in N_1  \mid t = \tcut (\lambda)\} \cup \FIX^{12}$. 

Consider further the case
  $t = \tcut (\lambda)$ for $(\lambda,t) \in N_1$. The following cases are possible: 
\begin{itemize}
\item $k \in (0, k_0)$. Then $\tcut(\lambda) = 4 K(k)$, thus $u_1 = \pi$, which implies   $x=0$. The equality $z=0$ holds only in the case $\cos u_2 = 0$. These points correspond to Maxwell points for the symmetry   $\varepsilon^2$, namely to the set   $\MAX^{20}$.
\item $k = k_0$. Then $\tcut(\lambda) = 4 K(k_0)$, i.e., $u_1 = u_{1z}(k_0)= \pi$, whence $x=z=0$. In this case $(\lambda,t) \in \MAX^{12}$, the points are Maxwell points for the both symmetries  $\varepsilon^1, \varepsilon^2$.
\item $k \in (k_0,1)$. Then  $\tcut(\lambda) = 2 p_z^1 (k)/ |\sigma|$, i.e., $u_1 = u_{1z} (k) \in (\pi/2, \pi)$, thus $z=0$. The equality $x=0$ holds only in the case $\sin u_2 = 0$. These points correspond to Maxwell points for the symmetry   $\varepsilon^1$, namely to the set $\MAX^{10}$.
\end{itemize}
Equality~(\ref{N00m}) follows.
\end{proofn}
 
The set $\FIX^{12}$ decomposes in two connected components:  $$\FIX_{\pm}^{12} = \{\nu \in N_7 \mid \cos \theta = \pm 1\}.$$

\begin{figure}[htbp]
\centering
\includegraphics[width=0.99\linewidth]{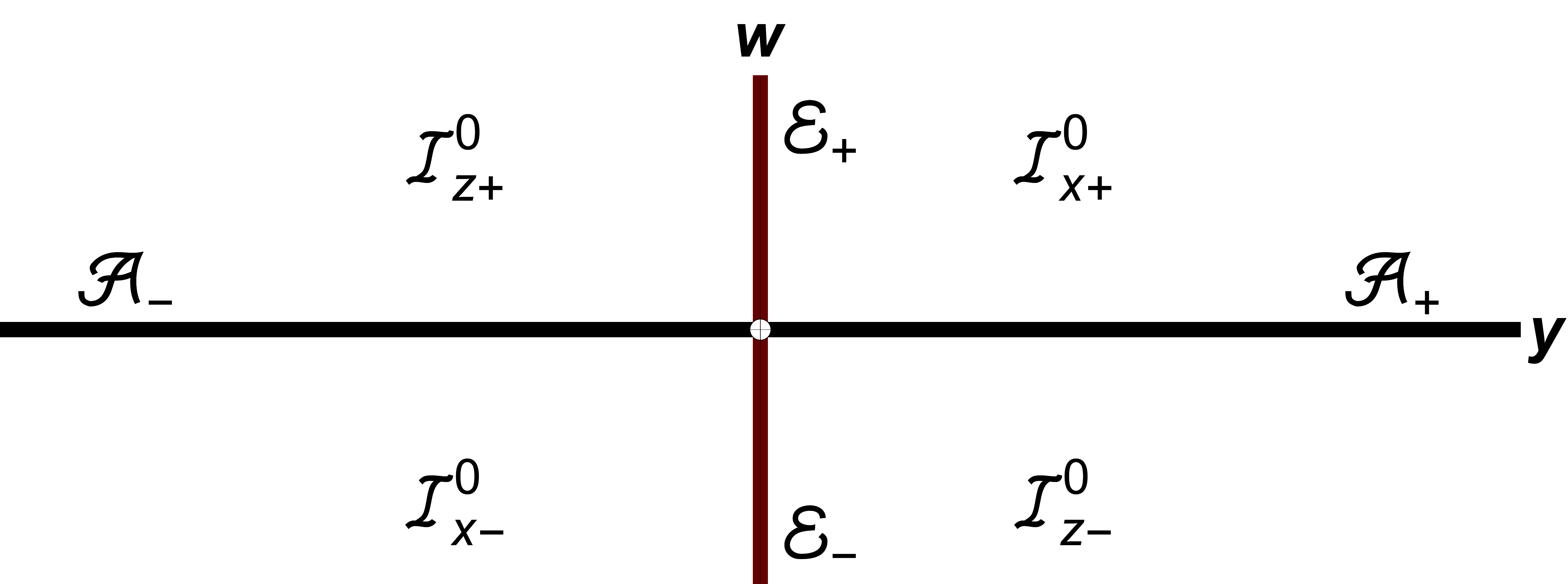}
\caption{Stratification of the plane   $\{x=z=0\}$.}\label{M00}
\end{figure}

A stratification of the set $M_{00}$ is shown in Table~\ref{tab:dec} and in Fig.~\ref{M00}.

\begin{table}[ht]
\caption{\label{tab:dec} Stratification of $M_{00}$ in coordinates $(y,w)$}
\begin{center}
\begin{tabular}{|c|c|c|c|c|c|c|c|c|}
\hline
$M_{0,0}$ & $\mathcal{I}_{x+}^0$ & $\mathcal{E}_+$ & $\mathcal{I}_{z+}^0$ & $\mathcal{A}_+$ & $\mathcal{A}_-$ & $\mathcal{I}_{z-}^0$ & $\mathcal{E}_-$ & $\mathcal{I}_{x-}^0$ \\
\hline
$y$ &  $(0,\infty)$ & 0 & $(-\infty,0)$ & $(0,\infty)$ & $(-\infty,0)$  & $(0,\infty)$  & 0 & $(-\infty,0)$ \\
\hline
$w$ &  \multicolumn{3}{c|}{$(0,+\infty)$} & \multicolumn{2}{c|}{$0$} & \multicolumn{3}{c|}{$(-\infty,0)$} \\
\hline
\end{tabular}
\end{center}
\end{table}

Lemmas~\ref{lemA+}, \ref{lemMAX12}  are obvious, and Lemmas~\ref{lemMAX20}, \ref{lemMAX10} are proved in Appendix A.

\begin{lem}\label{lemA+}
The mapping $\Exp \colon \FIX^{12}_{\pm} \to \mathcal{A}_{\pm}$ is a diffeomorphism.
\end{lem}

\begin{lem}\label{lemMAX12}
For any $u_2^0 \in [0, 2 \pi)$ the mapping $$\Exp \colon \big(\MAX_+^{12} \cap \{u_2 = u_2^0\}\big) \to \mathcal{E}_+$$ is a diffeomorphism.
\end{lem}

\begin{lem} \label{lemMAX20}
The mapping $\Exp \colon \MAX_{++}^{20} \to \mathcal{I}_{x+}^0$ is a diffeomorphism.
\end{lem}

\begin{lem} \label{lemMAX10}
The mapping $\Exp \colon \MAX_{++}^{10} \to \mathcal{I}_{z+}^0$ is a diffeomorphism.
\end{lem}

\begin{cor} \label{conj00}
The following restrictions of the exponential mapping are diffeomorphisms:
\begin{enumerate}
\item $\Exp \colon \FIX^{12}_- \to \mathcal{A}_-$,
\item $\Exp \colon \big(\MAX_-^{12} \cap \{u_2 = u_2^0\}\big) \to \mathcal{E}_-, u_2^0 \in [0, 2 \pi)$,
\item $\Exp \colon \MAX_{+-}^{20} \to \mathcal{I}_{x+}^0, \ \Exp \colon \MAX_{-+}^{20} \to \mathcal{I}_{x-}^0, \ \Exp \colon \MAX_{--}^{20} \to \mathcal{I}_{x-}^0,$
\item $\Exp \colon \MAX_{+-}^{10} \to \mathcal{I}_{z+}^0, \ \Exp \colon \MAX_{-+}^{10} \to \mathcal{I}_{z-}^0, \ \Exp \colon \MAX_{--}^{10} \to \mathcal{I}_{z-}^0.$
\end{enumerate}
\end{cor}
\begin{proofn}
Follows immediately from Lemmas~\ref{lemA+}--\ref{lemMAX10} respectively, with account of the symmetries   $\varepsilon^1, \varepsilon^2, \varepsilon^4$.
\end{proofn}

We sum up results of this section in the following statement.

\begin{teo}\label{thm:x=z=0}
There is a stratification
$$
\Cut \cap M_{00} = \{ q \in M \mid x=z=0, w\neq 0\} = \bigsqcup_{i\in \{+,-\}} \big( \mathcal{I}_{xi}^0 \sqcup \mathcal{I}_{zi}^0 \sqcup \mathcal{E}_i \big).
$$
To each point of the quadrants
 $\mathcal{I}_{x\pm}^0, \mathcal{I}_{z\pm}^0$ come exactly two sub-Riemannian minimizers, while to any point of the rays   $\mathcal{E}_\pm$ come a one-parametric family of minimizers.   

The rays $\mathcal{A}_+ \sqcup \mathcal{A}_- = M_{00} \backslash \Cut$ are filled by abnormal trajectories.  
\end{teo}

Now we describe minimizers for
  $q_1 = (0,y_1,0,w_1) \in M_{00}$: 
\begin{enumerate}
\item for $q_1 \in \mathcal{A}_{\pm}$ there is a unique minimizer   $\nu \in \FIX_{\pm}^{12}.$ 
\item for $q_1 \in \mathcal{I}_{x\pm}^0$ there are two minimizers 
\begin{align*}
\nu_1 &= (\hat{k}, \pi, \pi/2, \hat{\sigma}) \in \MAX_{\pm+}^{20}, \quad \nu_2 = (\hat{k}, \pi, 3\pi/2, \hat{\sigma}) \in \MAX_{\pm-}^{20}, \\
\hat{k} &\in (0, k_0), \quad \pm \hat{\sigma}>0.
\end{align*}
Figure~\ref{max1} shows two minimizers with   $\hat{k}=0.84, \hat{\sigma}=1$. 

\begin{figure}[htbp]
\centering
\includegraphics[width=0.99\linewidth]{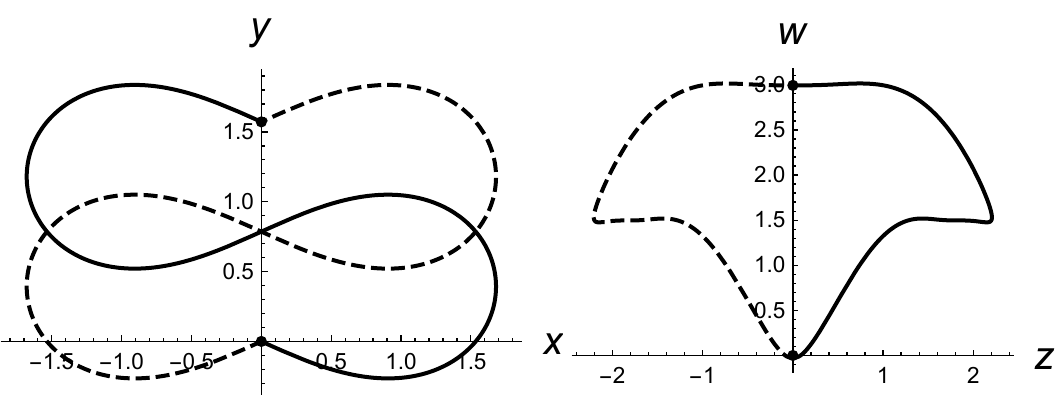}
\caption{Example of symmetric trajectories for   $x_1=z_1=0, y_1>0, w_1>0$. Left: projections to the plane   $(x,y)$, right: projections to the plane $(z,w)$.}\label{max1}
\includegraphics[width=0.99\linewidth]{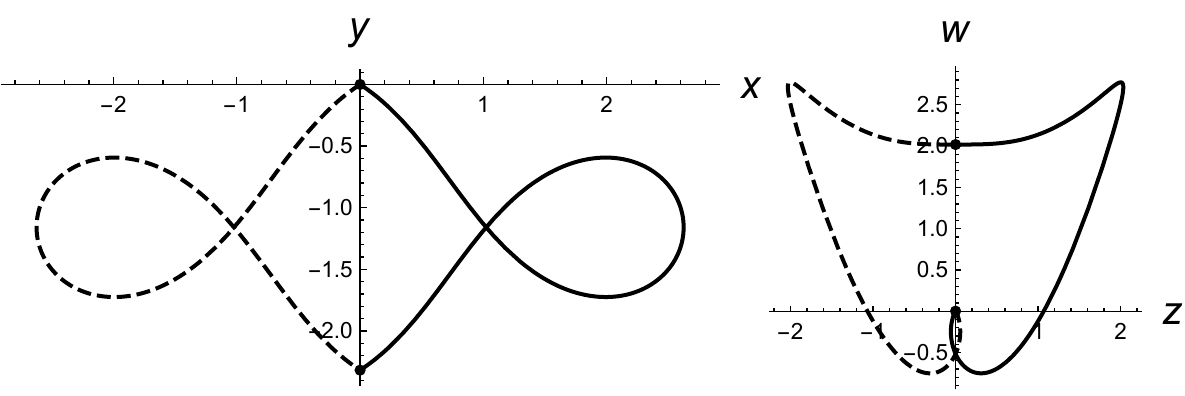}
\caption{Example of symmetric trajectories for   $x_1=z_1=0, y_1<0, w_1>0$. Left: projections to the plane   $(x,y)$, right: projections to the plane $(z,w)$.}\label{max2}
\includegraphics[width=0.85\linewidth]{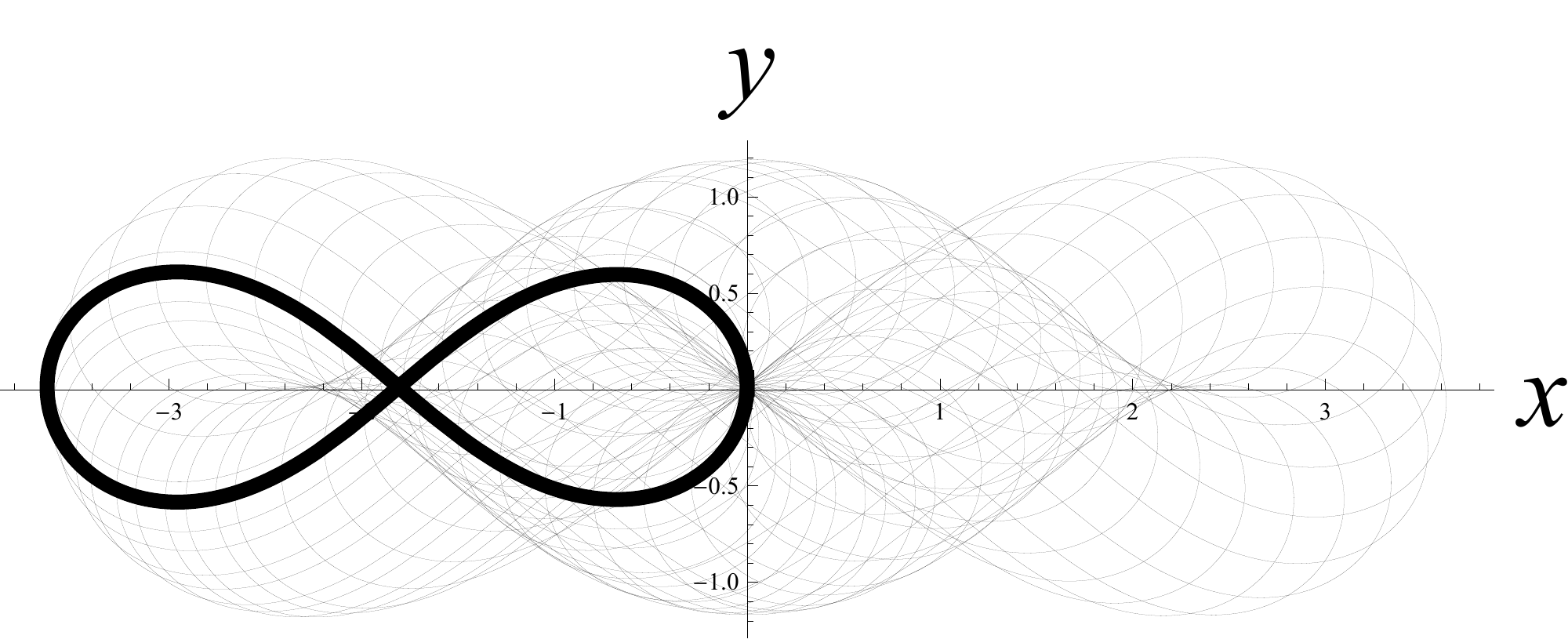}
\caption{One-parameter family of~"figure-of-eight" trajectories coming to the same point (projections to  $(x,y)$).}\label{eight1}
\end{figure}

\item for $q_1 \in \mathcal{I}_{z\pm}^0$ there are two minimizers   
\begin{align*}
\nu_1 &= \big(\hat{k}, u_{1z}(\hat{k}), 0, \hat{\sigma}\big) \in \MAX_{\pm+}^{10}, \quad \nu_2 = \big(\hat{k}, u_{1z}(\hat{k}), \pi, \hat{\sigma} \big) \in \MAX_{\pm-}^{10}, \\
\hat{k} &\in (k_0,1), \quad \pm \hat{\sigma}>0.
\end{align*}
Figure~\ref{max2} shows two optimal trajectories with   $\hat{k}=0.95, \hat{\sigma}=1$. 

\item for $q_1\in \mathcal{E}_\pm$ there exists a one-parameter family of minimizers  
\begin{align*}
\nu (u_2^0) &= (k_0, \pi, u_2^0, \hat{\sigma}) \in \MAX_{\pm}^{12}, \qquad u_2^0 \in [0, 2 \pi).
\end{align*}

In the plane $(x,y)$ the minimizers project to "figure-of-eight" elasticae   (see Fig.~\ref{eight1}, the boldface curve corresponds to   $u_2^0=0$). 

Despite the fact that solutions in the plane 
$(x,y)$ have the same form and are transformed one to another by parallel translations, there is no continuous symmetry in the space   $M$ that transforms these solutions one to another. There are only discrete symmetries that decompose solutions into four-tuples for   $u_2^0 \neq \pi n/2, n=0,\dots,3,$ and two pairs for    $u_2^0 = 0, \pi$ and for $u_2^0 = \pi/2, 3\pi/2$. Figure~\ref{eight2} shows a four-tuple of symmetric solutions. 

\begin{figure}[ht]
\includegraphics[width=0.99\linewidth]{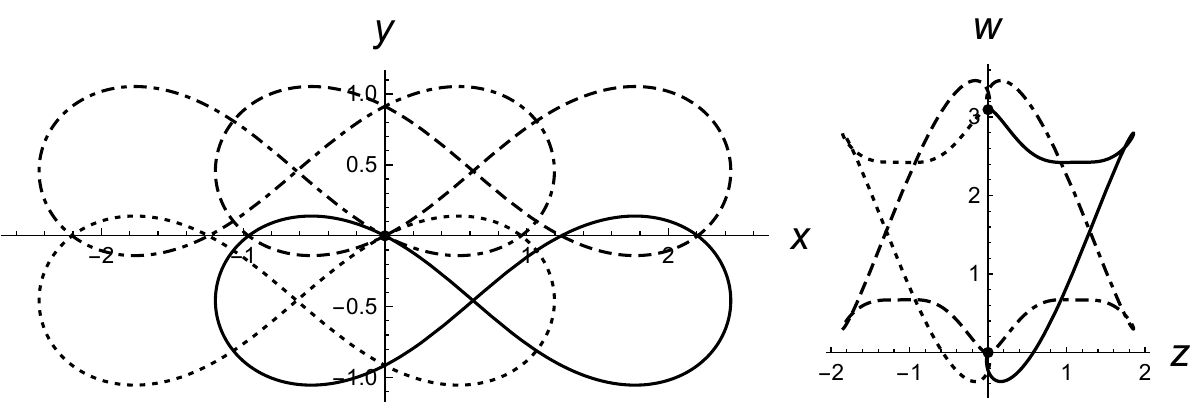}
\caption{Example of symmetric trajectories for   $x_1=y_1=z_1=0$. Left: projections to   $(x,y)$, right: projections to   $(z,w)$.}\label{eight2}
\end{figure}
\end{enumerate}

\section{Structure of $\Cut \cap \{z=0, \ x > 0\}$} \label{sec:z=0+}
In this section we study the set
  $\Cut \cap M_{+0}$. Consider the corresponding set in the preimage of the exponential mapping   $N_{+0} = \Exp^{-1}(M_{+0}) \cap \widehat{N}$ and study its intersection with the cut locus   $N_{\cut}$. Notice that the equality   $z=0$, by Lemma~\ref{lem0}, implies that    $N_{+0} \subset N_1 \cup N_7$. Since $\restr{\tcut}{C_7} = \infty$, we get
\begin{align*}
N_{+0} \cap N_{\cut} &= \{(k, u_1,u_2, \sigma)\in N_1 \mid k \in (k_0,1), u_1 = u_{1z}(k), \sin u_2 \neq 0 \} \\
&= \Big(\bigsqcup_{i,j \in \{+,-\}} \MAX_{ij}^{1+}\Big)  \sqcup \Big(\bigsqcup_{i \in \{+,-\}} \CMAX_i^{1+}\Big),
\end{align*}
where the sets $\MAX_{+\pm}^{1+}$,$\MAX_{-\pm}^{1+}$, and $\CMAX_{\pm}^{1+}$ are defined via Table~\ref{tab:z=0}.
\begin{table}[ht]
\caption{\label{tab:z=0} Components of $N_{+0} \cap N_{\cut} \cap N_1$} 
\begin{center}
\begin{tabular}{|c|c|c|c|c|c|c|c|c|c|c|c|}
\hline
& $\MAX_{++}^{1+}$ & $\CMAX_{+}^{1+}$ & $\MAX_{+-}^{1+}$  & $\MAX_{-+}^{1+}$ & $\CMAX_{-}^{1+}$ & $\MAX_{--}^{1+}$ \\
\hline
$u_2$ & $(0,\pi/2)$ & $\pi/2$ & $(\pi/2,\pi)$ & $(\pi, 3 \pi/2)$ & $3 \pi/2$ & $(3 \pi/2,2 \pi)$ \\
\hline
$\sigma$ & \multicolumn{3}{c|}{$(-\infty,0)$} & \multicolumn{3}{c|}{$(0,+\infty)$} \\
\hline
$k$ & \multicolumn{6}{c|}{$(k_0,1)$}  \\
\hline
$u_1$ & \multicolumn{6}{c|}{$u_{1z}(k)$} \\
\hline
\end{tabular}
\end{center}
\end{table} 

In order to describe a decomposition of the set
  $M_{+0}$, we pass to new coordinates  $\ds (Y^1,W^1) = 
\bigg(\frac{y}{x},\frac{w}{x^3}\bigg)$ invariant w.r.t. dilations $e^{t X_0}$ .
 
We use below the following notation:
\begin{align}
s_i &= \sin u_i, &&c_i = \cos u_i, &&d_i = \sqrt{1-k^2 s_i^2}, \label{scd} \\
E_1 &= E\big(u_1,k\big), &&F_1 = F\big(u_1,k\big), &&\Delta=1-k^2 s_1^2 s_2^2. \label{EFD}
\end{align}

In the case
  $N_{+0} \cap N_{\cut}$ the exponential mapping is determined by the formulas: 
\begin{align}
Y_1^1(k, u_1, u_2)& =- \frac{1 + k^2 (s_1^2 - 2) s_2^2}{2 k \, c_1 s_2 d_2}, \nonumber\\
W_1^1 (k, u_1, u_2)& = \frac{1}{48 k^3 s_1^3 c_1 d_1^3 s_2^3 d_2^3} \bigg(- E_1 c_1 \Delta^3 \nonumber\\
&\qquad + d_1^3 s_1 \Big(1  - k^2 s_1^2 s_2^2 \big(6 - 3 k^2 (4 - s_1^2) s_2^2 + 4 k^4 (2 - s_1^2) s_2^4\big)\Big)\bigg).\nonumber
\end{align}

The sets of conjugate points
  $\mathcal{CI}_{z\pm}^{+} = \Exp(\CMAX_{\pm}^{1+})$ consist of limit points for Maxwell points, these sets are parametrically defined as follows: 
\begin{align*}
\mathcal{CI}_{z\pm}^{+} = \{(Y^1, W^1) = \big(Y_1^1 (k, u_1, u_2), W_1^1 (k, u_1, u_2)\big) \mid \ & k\in (k_0,1), u_1 = u_{1z}(k),\\
												& \sin u_2 = \pm 1 \}. 
\end{align*}

Lemmas~\ref{diffzc}--\ref{diffzm} are proved in Appendix B.

\begin{lem}\label{diffzc} 
The mapping $\Exp\colon \CMAX_{+}^{1+} \to \mathcal{CI}_{z+}^{+}$ is a diffeomorphism.
\end{lem}

\begin{lem}\label{lem:c-z=0+} 
\begin{itemize}
\item[$1)$] The curve $\mathcal{CI}_{z+}^{+}$ in the plane $(Y^1,W^1)$ is a graph of a smooth function   $W_{\conj}^{1}(Y^1)$, increasing from   $-\infty$ to $\infty$ in the interval   $Y^1 \in (Y_0^1, \infty)$, where $\ds Y_0^1 = \frac{1 - 2 k_0^2}{2 k_0 \sqrt{1 - k_0^2}}<0$. 
\item[$2)$] The curve $W_1 = W_{\conj}^{1}(Y^1)$ lies below the line  $W^1 = Y_1/6$.
\item[$3)$] $\ds \lim\limits_{Y^1\to +\infty} \frac{W_{\conj}^{1}(Y^1)}{Y^1} \in \R \backslash \{0\}$.
\end{itemize}
\end{lem}

Define the following sets in the plane $(Y^1,W^1)$: 
\begin{align*}
\mathcal{I}_{z+}^+ &= \{(Y^1,W^1) \mid  W_{\conj}^{1}(Y^1) > W^1, Y^1 \in (Y_0^1,\infty)\}, \\
\mathcal{I}_{z-}^+ &= \{(Y^1,W^1) \in \R^2 \mid (-Y^1,-W^1) \in \mathcal{I}_{z+}^+\},
\end{align*}
and show their relation to the cut locus.

\begin{lem}\label{diffzm}
The mapping $\Exp\colon \MAX_{+-}^{1+} \to \mathcal{I}_{z+}^+$ is a diffeomorphism.
\end{lem}

\begin{cor} \label{conject01}
The following restrictions of the exponential mapping are diffeomorphisms:
$$\Exp \colon \MAX_{++}^{1+} \to \mathcal{I}_{z+}^+, \quad	  \Exp\colon \CMAX_{-}^{1+} \to \mathcal{CI}_{z-}^{+}, \quad  \Exp \colon \MAX_{-\pm}^{1+} \to \mathcal{I}_{z-}^+.$$
\end{cor}
\begin{proofn}
Follows immediately from Lemmas~\ref{diffzc},~\ref{diffzm}, with account of the symmetries  $\varepsilon^1$, $\varepsilon^2$, $\varepsilon^4$.
\end{proofn}

\begin{figure}[ht]
\centering
\includegraphics[width=0.99\linewidth]{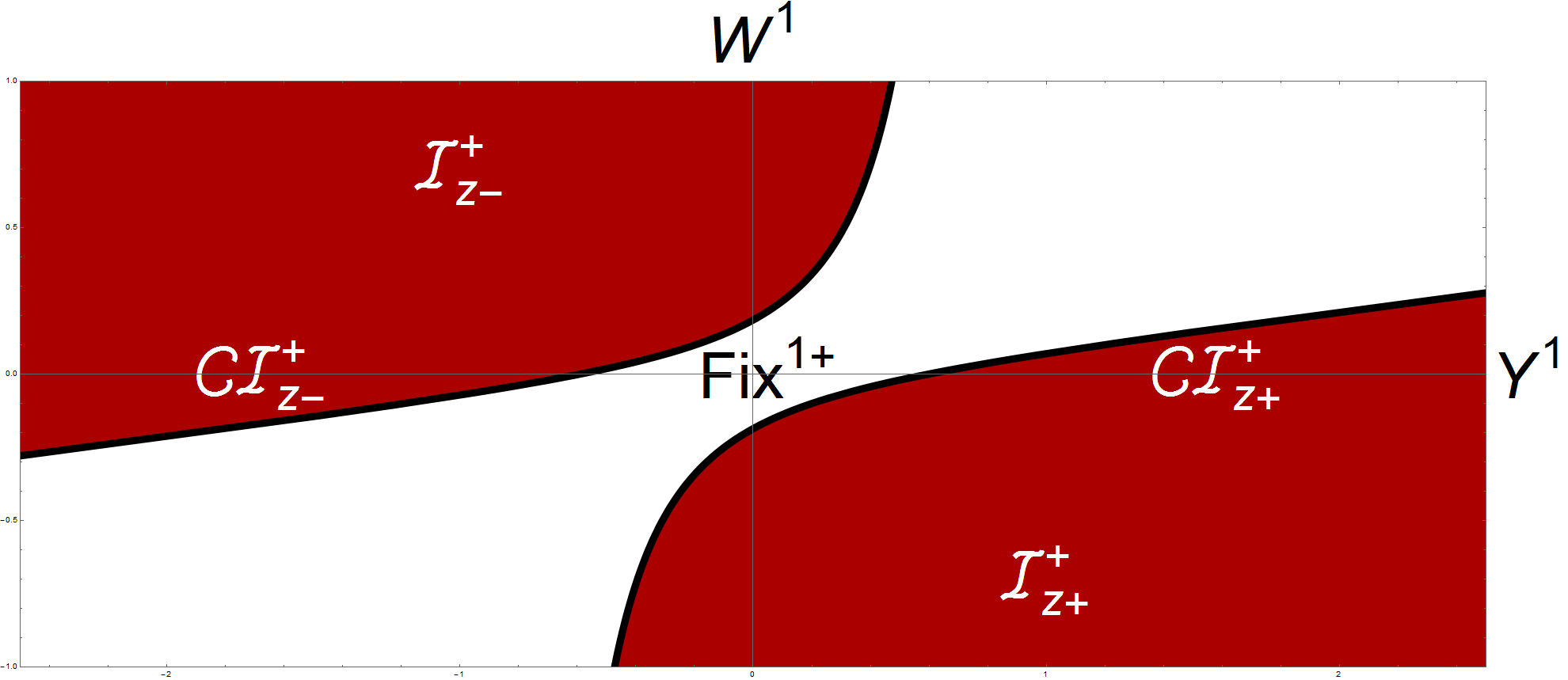}
\caption{Decomposition of the set $M_{+0}$}\label{decom1y0w}
\end{figure}

Figure~\ref{decom1y0w} shows a decomposition of the plane   $(Y^1,W^1)$.

\begin{teo}\label{thm:z=0+}
There is a stratification
\begin{align*}
\Cut \cap M_{+0} &= \mathcal{I}_{z+}^+ \sqcup \mathcal{CI}_{z+}^{+}  \sqcup \mathcal{I}_{z-}^+ \sqcup \mathcal{CI}_{z-}^{+}.
\end{align*}
Moreover,
\begin{align}
\mathcal{I}_{z+}^+ &= \Big\{q \in M \mid z=0, \  x>0, \ y > Y_0^1 x, \ w < W_{\conj}^1 (y/x) x^3\Big\} \cong \R^3, \label{I+z+} \\
\mathcal{I}_{z-}^+  &= \Big\{q \in M \mid z=0, \ x>0, \ y < - Y_0^1 x, \ w > - W_{\conj}^1 (-y/x) x^3\Big\} \cong \R^3, \nonumber \\
\mathcal{CI}_{z+}^{+} &= \Big\{q \in M \mid z=0, \ x>0, \ y > Y_0^1 x, \ w = W_{\conj}^1 (y/x) x^3\Big\} \cong \R^2, \nonumber \\
\mathcal{CI}_{z-}^{+} &= \Big\{q \in M \mid z=0, \ x>0, \ y < - Y_0^1 x, \ w = - W_{\conj}^1 (-y/x) x^3\Big\} \cong \R^2. \nonumber 
\end{align}
For each point of the sets
  $\mathcal{I}_{z\pm}^+$ there exist two minimizers, and for each point of the remaining part   $M_{+0}\backslash(\mathcal{I}_{z+}^+ \sqcup \mathcal{I}_{z-}^+)$ there is a unique minimizer.  
\end{teo}

\begin{proofn}
Follows from Lemmas~\ref{diffzc}--\ref{diffzm}, and from Corollary~\ref{conject01}.
\end{proofn}


\section{Structure of $\Cut \cap \{z=0\}$} \label{sec:z=0}
In this section we describe intersection of the cut locus with the subspace 
 $M_z = \{q\in M \mid z =0 \}$.

\begin{teo}
There is a stratification
\begin{equation}\label{CutMz}
\Cut \cap M_z = \bigsqcup_{i\in \{+,-\}} \Big(\mathcal{I}_{zi} \sqcup \big(\sqcup_{j \in \{+,-\}} \mathcal{CI}_{zi}^j\big) \sqcup \mathcal{I}_{xi}^0 \sqcup \mathcal{E}_i\Big),
\end{equation}
where $\mathcal{I}_{zi} = \bigsqcup_{j\in \{+,-,0\}} \mathcal{I}_{zi}^j$, \ $\mathcal{I}_{zi}^- =\varepsilon^2(\mathcal{I}_{zi}^+)$, \ $\mathcal{CI}_{zi}^- = \varepsilon^2 (\mathcal{CI}_{zi}^+), \ i \in \{+,-\}$. Moreover,
\begin{align}
\mathcal{I}_{z+} &= \{q \in M \mid z=0, \ y>Y_0^1 |x|, \ w <G_1 (x,y)\} \cong \R^3,  \label{Iz+}\\
\mathcal{I}_{z-} &= \varepsilon^4(\mathcal{I}_{z+})\cong \R^3, \nonumber
\end{align} 
where  $G_1$ is a function continuous on the set   $\{(x,y)\in \R^2 \mid y > Y_0^1 |x|\}$  and satisfying the properties
\begin{equation*}
G_1(0,y) = 0, \quad G_1(-x,y) = G_1(x,y), \quad  G_1(\rho x, \rho y) = \rho^3 G_1(x, y), \qquad \rho>0.
\end{equation*}
\end{teo} 
\begin{proofn}
By virtue of the equalities $\varepsilon^2(\Cut) = \Cut, \varepsilon^2(M_{+0}) = M_{-0},$  we get, with account of Theorem~\ref{thm:z=0+}, a decomposition
$$
\Cut \cap M_{-0} = \varepsilon^2 (\Cut \cap M_{+0}) = \bigsqcup_{i\in\{+,-\}} (\mathcal{I}_{zi}^- \sqcup \mathcal{CI}_{zi}^-).$$
Whence, with account of Theorem~\ref{thm:x=z=0}, we get a stratification~(\ref{CutMz}).

Representation~(\ref{Iz+}) is obtained from~(\ref{I+z+}), Table~\ref{tab:dec}, and the equalities
$$ 
\mathcal{I}_{z+}^-= \{q\in M \mid z=0, \ x<0, \ y>Y_0^1 |x|, \ w<W_{\conj}^1 (y/|x|) |x|^3\}
$$
for the function 
$$
G_1 (x,y) = \left\{\begin{array}{l}
			W_{\conj}^1 (y/|x|) |x^3|,  \  x\neq0, y>Y_0^1 |x|, \\
			0, \ x=0, y>0.
			\end{array}\right.
$$
Continuity of the function
  $G_1$ on the set $\{(x,y) \mid x\neq 0, y> Y_0^1 |x|\}$ follows from continuity of the function  $W_{\conj}^1$ on the ray $(Y_0^1, +\infty)$, see Lemma~\ref{lem:c-z=0+}. In order to prove continuity of the function   $G_1$ on the ray $\{(x,y) \mid x=0, y>0\}$, take any sequence   $(x_n,y_n)$, $x_n \to +0$, $y_n \to \bar{y} >0$. Then, with account of item  3) of Lemma~\ref{lem:c-z=0+}, we obtain
$$
G_1(x_n,y_n) = W_{\conj}^1(y_n/x_n) x_n^3 = \frac{W_{\conj}^1(y_n/x_n)}{y_n/x_n} y_n x_n^2 \to 0.
$$
Thus the function 
  $G_1$ is continuous on its whole domain   $\{(x,y) \mid y>Y_0^1 |x|\}$.

Representation~(\ref{Iz+}) implies that the stratum   $\mathcal{I}_{z+}$ is homeomorphic to   $\R^3$. The theorem is proved.
\end{proofn}

\section{Structure of $\Cut \cap \{x=0, \ z>0 \}$}\label{sec:x=0+}
In this section we describe the intersection
  $M_{0+} \cap \Cut$. We study the set $N_{0+}=\Exp^{-1}(M_{0+})\cap \widehat{N}$. It follows from Theorem~\ref{tcut} that
$$N_{0+}\cap N_{\cut} \subset N_1 \cup N_2 \cup N_6.$$

Now we consider cut points for each of the subsets
  $N_i, i =1,2,6$. 

\subsection{Subcase $N_1$}
We have
\begin{align*}
N_{0+}\cap N_{\cut} \cap N_1 &= \{(k,u_1,u_2,\sigma)\in N_1 \mid k\in(0,k_0), u_1=\pi\} \\
&= \Big(\bigsqcup_{i,j \in \{+,-\}} \MAX_{1ij}^{2+}\Big) \sqcup \Big(\bigsqcup_{i \in \{+,-\}} \CMAX_{1i}^{2+}\Big),
\end{align*}
where the sets $\CMAX_{1\pm}^{2+}, \MAX_{1+\pm}^{2+}, \MAX_{1-\pm}^{2+}$ are defined by values of the parameters   $u_2, \sigma$ via Table~\ref{tab:C1-x=0}; notice that  inequality   $z>0$ implies $\cos u_2 > 0$.

\begin{table}[ht]
\caption{\label{tab:C1-x=0} Components of $N_{0+} \cap N_{\cut} \cap N_1$}
\begin{center}
\begin{tabular}{|c|c|c|c|c|c|c|c|c|c|c|c|}
\hline
 & $\CMAX_{1+}^{2+}$ & $\MAX_{1++}^{2+}$ & $\MAX_{1+-}^{2+}$ & $\CMAX_{1-}^{2+}$ & $\MAX_{1-+}^{2+}$ & $\MAX_{1--}^{2+}$ \\
\hline
$u_2$ & 0 & $(0,\pi/2)$ & $(3\pi/2,2 \pi)$ & $0$ & $(0,\pi/2)$ & $(3\pi/2,2 \pi)$ \\
\hline
$\sigma$ & \multicolumn{3}{c|}{$(0, \infty)$} & \multicolumn{3}{c|}{$(-\infty,0)$} \\
\hline
$k$ & \multicolumn{6}{c|}{$(0, k_0)$} \\
\hline
$u_1$ & \multicolumn{6}{c|}{$\pi$} \\
\hline
\end{tabular}
\end{center}
\end{table} 
In order to define a decomposition in the set
 $M_{0+}$, we introduce new coordinates
  $\ds (Y^2,W^2) = \bigg(\frac{y}{\sqrt{z}},\frac{w}{\sqrt{z^3}}\bigg)$ invariant w.r.t. dilations $e^{t X_0}$.  Then the exponential mapping takes the form:
\begin{align*}
Y_1^2(k,u_2) = \sqrt{\frac{2 \iota_1(k)}{k c_2}},  \qquad  W_1^2(k,u_2) = \frac{\iota_2(k) + k^2 \iota_1(k) (1 + 3 c_2^2)}{3 \big(2 k \iota_1(k) c_2\big)^{3/2}},
\end{align*}
where the functions $\iota_1(k), \iota_2(k)$ are defined in Appendix A, see~(\ref{iota12}).  

The sets of conjugate points
  $\mathcal{CI}_{x\pm}^{+} = \Exp(\CMAX_{1\pm}^{2+})$ are defined as follows:  
\begin{align}
\mathcal{CI}_{x\pm}^{+} = \{(Y^2,W^2) = (\pm Y_1^2(k,u_2), \pm W_1^2(k,u_2)) \mid k \in (0, k_0), u_2 = 0\}.
\end{align} 

Lemmas~\ref{lem:c2-x=0}--\ref{diffxm2} and \ref{lem:smoothx2} are proved in Appendix C.
Lemma~\ref{lem:Max62} is obvious.

\begin{lem}\label{lem:c2-x=0}
The mapping $\Exp: \CMAX_{1+}^{2+} \to \mathcal{CI}_{x+}^{+}$ is a diffeomorphism.
\end{lem}

\begin{lem}\label{lemcon21}
The curve $\mathcal{CI}_{x+}^{+}$ in the plane $(Y^2,W^2)$ is a graph of a certain smooth function   $W_{\conj}^{21}(Y^2)>0$, decreasing from $\infty$ to $0$ at the interval $Y^2\in (0,\infty)$.
\end{lem}

Define the following sets in the plane $(Y^2,W^2)$:
\begin{align}
\mathcal{I}_{x+}^{+} &= \{(Y^2,W^2) \mid W_{\conj}^{21}(Y^2)<W^2, \ Y^2 \in (0, \infty)\}, \label{Ix++}\\
\mathcal{I}_{x-}^+ &= \{(Y^2,W^2) \in \R^2 \mid (-Y^2,-W^2) \in \mathcal{I}_{x+}^+\} \nonumber
\end{align}
and show their relation to the cut locus.

\begin{lem} \label{diffxm1}
The mapping  $\Exp \colon \MAX_{1++}^{2+} \to \mathcal{I}_{x+}^+$ is a diffeomorphism.
\end{lem}


\begin{cor} \label{conj02}
The following restrictions of the exponential mapping are diffeomorphisms:
$$\Exp \colon \MAX_{1+-}^{2+} \to \mathcal{I}_{x+}^+, \quad \Exp\colon \CMAX_{1-}^{2+} \to \mathcal{CI}_{x-}^{+}, \quad \Exp \colon \MAX_{1-\pm}^{2+} \to \mathcal{I}_{x-}^+.$$
\end{cor}
\begin{proofn}
Follows immediately from Lemmas~\ref{lem:c2-x=0}, \ref{diffxm1}, with account of the symmetries    $\varepsilon^1, \varepsilon^2, \varepsilon^4$.
\end{proofn}

\subsection{Subcase $N_2$.}
There holds a representation
\begin{align*}
N_{0+}\cap N_{\cut} \cap N_2 &= \{(k,u_1,u_2,\sigma)\in N_2 \mid k\in(0,1), u_1=\pi/2\}\\
   &=\bigsqcup_{i,j \in \{+,-\}} \big(\CMAX_{2ij}^{2+} \sqcup \MAX_{2ij}^{2+}\big),
\end{align*}
where $\CMAX_{2+\pm}^{2+}, \CMAX_{2-\pm}^{2+}, \MAX_{2+\pm}^{2+}$ and $\MAX_{2-\pm}^{2+}$ are defined by values of the parameters   $u_2, \sigma$ via Table~\ref{tab:C2-x=0}. Notice that the inequality   $z>0$ implies that $\sgn c = 1$.

\begin{table}[ht]
\caption{\label{tab:C2-x=0} Components of $N_{0+} \cap N_{\cut} \cap N_2$}
\begin{center}
\begin{tabular}{|c|c|c|c|c|c|c|c|c|c|c|c|c|c|}
\hline
 & \begin{sideways}$\CMAX_{2++}^{2+}$\end{sideways} & \begin{sideways}$\MAX_{2++}^{2+}$\end{sideways} & \begin{sideways}$\CMAX_{2+-}^{2+}$\end{sideways} &  \begin{sideways}$\MAX_{2+-}^{2+}$\end{sideways} &  \begin{sideways}$\CMAX_{2-+}^{2+}$\end{sideways} & \begin{sideways}$\MAX_{2-+}^{2+}$\end{sideways} & \begin{sideways}$\CMAX_{2--}^{2+}$\end{sideways} &  \begin{sideways}$\MAX_{2--}^{2+}$\end{sideways}  \\
\hline
$u_2$ & $0$ & $(0,\pi/2)$ & $\pi/2$ &  $(\pi/2,\pi)$ & $0$ & $(0,\pi/2)$ & $\pi/2$ & $(\pi/2,\pi)$ \\
\hline
$\sigma$ & \multicolumn{4}{c|}{$(0, \infty)$} & \multicolumn{4}{c|}{$(-\infty,0)$} \\
\hline
$k$ & \multicolumn{8}{c|}{$(0, 1)$} \\
\hline
$u_1$ & \multicolumn{8}{c|}{$(0, \pi/2)$} \\
\hline
\end{tabular}
\end{center}
\end{table}

For $\nu \in N_{0+} \cap N_{\cut} \cap N_2$ the exponential mapping is defined in the coordinates   $(Y^2,W^2)$ by the formulas:
\begin{align*}
Y_2^2(k,u_2) &= - \sqrt{\frac{\iota_4(k) d_2}{\sqrt{1 - k^2}}} < 0,  \\
W_2^2(k,u_2) &= \frac{k^4 K(k) d_2^2 - \iota_4(k) \big(8 - 7 k^2 - k^2 (2 - k^2) s_2^2\big)}{12 \sqrt{\iota_4^3(k) (1 - k^2)^{3/2} d_2}},
\end{align*}
where the function $\iota_4(k)$ is defined in Appendix C, see~(\ref{iota4}).  

The sets of conjugate points
  $\mathcal{CN}_{x\pm+}^{+}$, $\mathcal{CN}_{x\pm-}^{+}$ are defined as follows:  
\begin{align}
\mathcal{CN}_{x\pm+}^{+} &= \Big\{(Y^2,W^2)=\big(\pm Y_2^2(k, 0), \pm W_2^2(k, 0)\big) \mid k \in (0, 1)\Big\}, \\
\mathcal{CN}_{x\pm-}^{+} &= \Big\{(Y^2,W^2)=\big(\pm Y_2^2(k, \pi/2), \pm W_2^2(k, \pi/2)\big) \mid k \in (0, 1)\Big\}.
\end{align} 

Lemmas~\ref{diffxc2+}-\ref{diffxm2} are proved in Appendix C.
\begin{lem}\label{diffxc2+}
The mapping $\Exp: \CMAX_{2++}^{2+} \to \mathcal{CN}_{x++}^{+}$ is a diffeomorphism.
\end{lem}

\begin{lem}\label{lemcon22+}
The curve $\mathcal{CN}_{x++}^{+}$ in the plane $(Y^2,W^2)$ is a graph of a smooth function   $W_{\conj}^{22+}(Y^2)$ decreasing from   $\infty$ to $1/\sqrt{\pi}$ at the interval $Y^2 \in  (-\infty, 0)$.
\end{lem}

\begin{lem}\label{diffxc2-}
The mapping $\Exp: \CMAX_{2+-}^{2+} \to \mathcal{CN}_{x+-}^{+}$ is a diffeomorphism.
\end{lem}

\begin{lem}\label{lemcon22-}
The curve $\mathcal{CN}_{x+-}^{+}$ in the plane $(Y^2,W^2)$ is a graph of a smooth function   $W_{\conj}^{22-}(Y^2)$, decreasing from   $0$ to $-1/\sqrt{\pi}$ at the interval $(-\infty, 0)$.
\end{lem}

\begin{cor} \label{con22-}
The following restrictions of the exponential mapping are diffeomorphisms:
$$
\Exp \colon \CMAX_{2-+}^{2+} \to \mathcal{CN}_{x-+}^{+}, \qquad
\Exp \colon \CMAX_{2--}^{2+} \to \mathcal{CN}_{x--}^{+}.
$$

The curves $\mathcal{CN}_{x-+}^{+}$, $\mathcal{CN}_{x--}^{+}$ in the plane $(Y^2,W^2)$ are respectively graphs of the functions $- W_{\conj}^{22+}(-Y^2)$, $- W_{\conj}^{22-}(-Y^2)$, $Y \in (0, \infty)$.
\end{cor}

Further we study the mutual disposition of the curves 
  $\mathcal{CI}_{x\pm}^{+}$, $\mathcal{CN}_{x\pm+}^{+}$, $\mathcal{CN}_{x\pm-}^{+}$. By virtue of Lemma~\ref{lemcon22+} and the symmetry $\varepsilon^4$~(\ref{vareps4}) it follows that the curves   $\mathcal{CN}_{x++}^{+}$ and $\mathcal{CN}_{x-+}^{+}$ belong respectively to the second and fourth quadrants of the plane   $(Y^2,W^2)$.
 
Lemma~\ref{lemcon21} and Corollary~\ref{con22-} imply that the curves   $\mathcal{CI}_{x+}^{+}, \mathcal{CN}_{x--}^{+}$ belong to the first quadrant. We show in the following lemma that they do not intersect each another.  

\begin{lem}\label{conjcross}
There holds inequality $-W_{\conj}^{22-}(- Y^2)<W_{\conj}^{21}(Y^2)$ for $Y^2 > 0$.
\end{lem} 

Define the following sets in the plane $(Y^2,W^2)$:
\begin{align*}
\mathcal{N}_{x+}^+ &= \{(Y^2,W^2)\in \R^2 \mid W_{\conj}^{22-}(Y^2)< W^2 < W_{\conj}^{22+}(Y^2), \ Y^2 < 0\}, \\
\mathcal{N}_{x-}^+ &= \{(Y^2,W^2)\in \R^2 \mid (-Y^2,-W^2) \in \mathcal{N}_{x+}^+\}.
\end{align*}

\begin{lem}\label{diffxm2}
The mapping $\Exp \colon \MAX_{2--}^{2+} \to \mathcal{N}_{x-}^+$ is a diffeomorphism.
\end{lem}

\begin{cor}
The following restrictions of the exponential mapping are diffeomorphisms:
$$ \Exp \colon \MAX_{2-+}^{2+} \to \mathcal{N}_{x-}^+, \  \Exp \colon \MAX_{2++}^{2+} \to \mathcal{N}_{x+}^+, \  \Exp \colon \MAX_{2+-}^{2+} \to \mathcal{N}_{x+}^+.$$
\end{cor}

\subsection{Subcase $C_6$}
There holds a decomposition
\begin{align*}
N_{0+}\cap N_{\cut} \cap N_6 &= \{(c, \theta, t) \in N_6 \mid t= 2 \pi/|c|, \theta \in S^1, c>0\} \\
&=\bigsqcup_{i\in \{+,-\}} \Big(\CMAX_{6i}^{2+} \sqcup \MAX_{6i}^{2+} \Big),  
\end{align*}
where the sets $\CMAX_{6\pm}^{2+}, \MAX_{6\pm}^{2+}$ are defined by values of the parameter   $\theta$ via Table~\ref{tabC6}.

\begin{table}[ht]
\caption{\label{tabC6} Components of $N_{0+} \cap \N_{\cut} \cap N_6$}
\begin{center}
\begin{tabular}{|c|c|c|c|c|}
\hline
 & $\MAX_{6+}^{2+}$ & $\CMAX_{6-}^{2+}$ & $\MAX_{6-}^{2+}$ &  $\CMAX_{6+}^{2+}$  \\
\hline
$\theta$ & $(-\pi,0)$ & $0$ & $(0,\pi)$ &  $\pi$\\
\hline
$c$ &  \multicolumn{4}{c|}{$(0, +\infty)$}\\
\hline
$t$ & \multicolumn{4}{c|}{$ 2\pi/|c|$}\\
\hline
\end{tabular}
\end{center}
\end{table} 

For $\nu \in N_{0+}\cap N_{\cut} \cap N_6$ the exponential mapping is defined in the coordinates   $(Y^2,W^2)$ by the formulas: 
$Y_6^2(\theta) = 0$,  $W_6^2(\theta) = -{\cos \theta}/{\sqrt{\pi}}$.

In the image of the exponential mapping the sets of conjugate points 
  $\mathcal{CC}_{x\pm}^{+}$ are defined as follows:
\begin{align*}
\mathcal{CC}_{x+}^{+} &= \{(Y^2,W^2)=\big(Y_6^2(\pi), W_6^2(\pi)\big) = (0,1/\sqrt{\pi})\}, \\
\mathcal{CC}_{x-}^{+} &= \{(Y^2,W^2)=\big(Y_6^2(0), W_6^2(0)\big) = (0,-1/\sqrt{\pi})\}.
\end{align*}

Define the following set in the plane
  $(Y^2,W^2)$:
\begin{align*}
\mathcal{C}_{x}^+ = \{(Y^2,W^2)\in \R^2 \mid Y^2 = 0, |W^2|< 1/\sqrt{\pi}\}.
\end{align*}

\begin{lem}\label{lem:Max62}
The mappings $\Exp (\MAX_{6\pm}^{2+}) \to \mathcal{C}_{x}^+$ are diffeomorphisms.
\end{lem}

\begin{figure}[ht]
\centering
\includegraphics[width=0.99\linewidth]{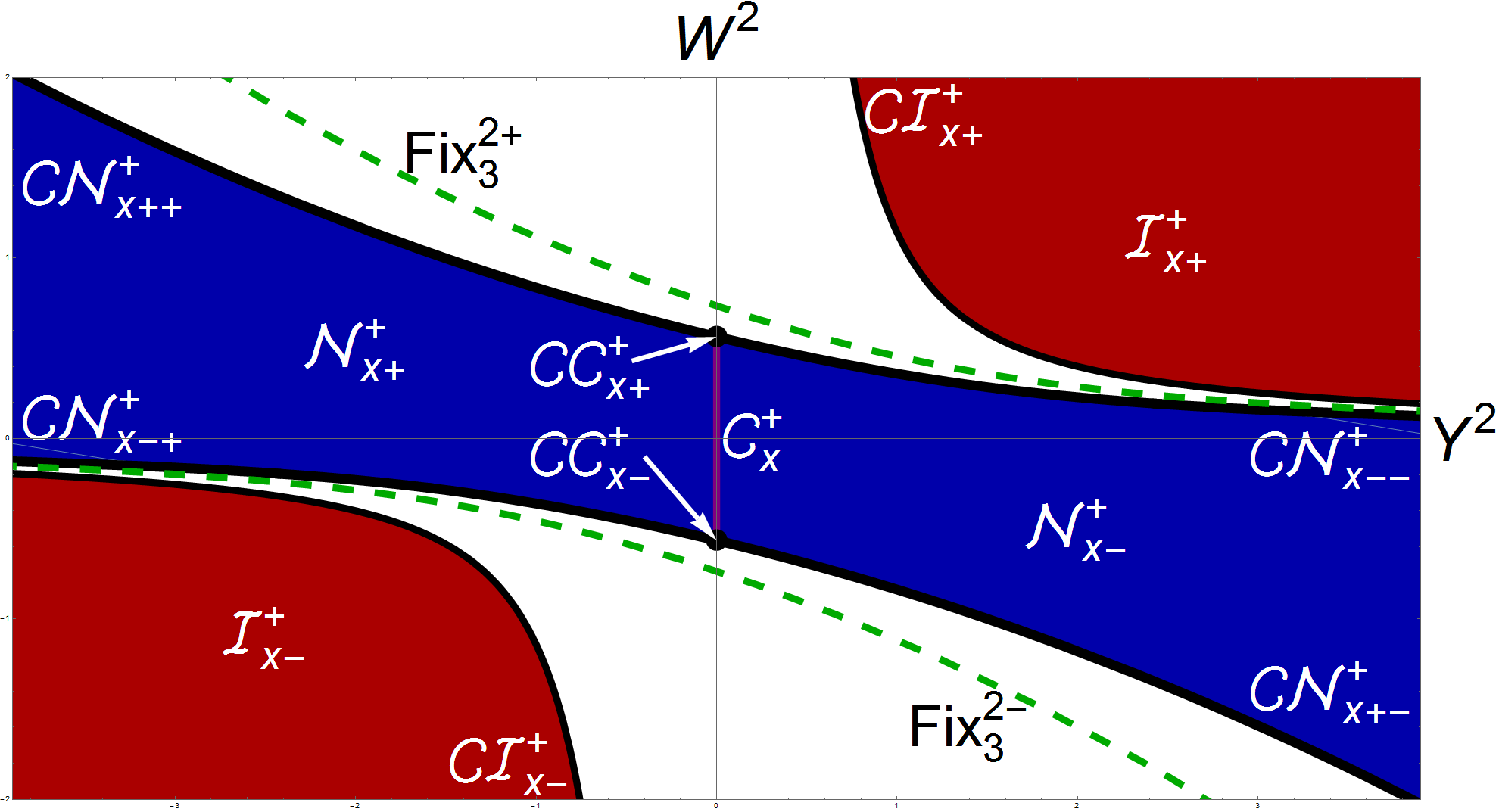}
\caption{Decomposition of the set   $M_{0+}$ }\label{fig:decx}
\end{figure}

A decomposition of the plane
  $(Y^2,W^2)$ is shown in Fig.~\ref{fig:decx}.

Denote the sets
\begin{align*}
\mathcal{N}_{x}^+ &= \mathcal{N}_{x+}^+ \sqcup \mathcal{C}_{x}^+ \sqcup \mathcal{N}_{x-}^+, \\
\mathcal{CN}_{x+}^+ &= \mathcal{CN}_{x++}^+ \sqcup \mathcal{CC}_{x+}^+ \sqcup \mathcal{CN}_{x--}^+,\qquad \mathcal{CN}_{x-}^+ = \mathcal{CN}_{x+-}^+ \sqcup \mathcal{CC}_{x-}^+ \sqcup \mathcal{CN}_{x-+}^+.
\end{align*}

\begin{lem} \label{lem:smoothx2}
The set $\mathcal{CN}_{x+}$ (resp. $\mathcal{CN}_{x-}^+$) forms a smooth curve in the upper (lower) half-plane, which is a graph of a smooth function $W_{\conj}^{22}(Y^2)$ (resp. $-W_{\conj}^{22}(-Y^2)$), where
$$W_{\conj}^{22} (Y^2) = \left\{\begin{array}{l}
			W_{\conj}^{22+} (Y^2),  \qquad   Y^2\in (-\infty,0), \\
			-W_{\conj}^{22-} (-Y^2), \   Y^2\in (\infty,0), \\
			1/\sqrt{\pi}, \qquad \qquad Y^2 = 0.
			\end{array}\right.$$
\end{lem}

We sum up results of this section as follows.

\begin{teo} \label{thm:x=0+}
There is a stratification
\begin{align*}
\Cut \cap M_{0+} =& \mathcal{I}_{x+}^+ \sqcup \mathcal{I}_{x-}^+ \sqcup \mathcal{N}_{x}^+ \sqcup \mathcal{CI}_{x+}^{+} \sqcup \mathcal{CI}_{x-}^{+} \sqcup \mathcal{CN}_{x+}^{+}\sqcup \mathcal{CN}_{x-}^{+}.
\end{align*}
Moreover, 
\begin{align}
\mathcal{I}_{x+}^{+} &= \{q\in M \mid x=0, \ z>0, \ y>0, \ w > W_{\conj}^{21}(y/\sqrt{z}) \sqrt{z^3}\} \cong \R^3, \label{I+x+}\\
\mathcal{I}_{x-}^{+} &= \{q\in M \mid x=0, \ z>0, \ y<0, \ w < -W_{\conj}^{21}(-y/\sqrt{z}) \sqrt{z^3}\} \cong \R^3,  \nonumber\\
\mathcal{N}_{x}^+ &= \{q \in M \mid x=0, \ z>0, \ - W_{\conj}^{22}(-y/\sqrt{z}) \sqrt{z^3} < w < W_{\conj}^{22}(y/\sqrt{z}) \sqrt{z^3}\} \cong \R^3, \nonumber \\
\mathcal{CI}_{x\pm}^+ &= \{q \in M \mid x=0, \ z>0, \ \pm y > 0, \ w = \pm W_{\conj}^{21}(\pm y/\sqrt{z}) \sqrt{z^3}\} \cong \R^2, \nonumber \\
\mathcal{CN}_{x\pm}^+ &= \{q \in M \mid x=0, \ z>0, \ w = \pm W_{\conj}^{22}(\pm y/\sqrt{z}) \sqrt{z^3} \} \cong \R^2.\nonumber 
\end{align}
To each point of the strata
  $\mathcal{I}_{x\pm}^{+}$, $\mathcal{N}_{x}^+$ come two minimizers, and to each point of the remaining part   $M_{0+} \slash (\mathcal{I}_{x+}^{+} \sqcup \mathcal{I}_{x-}^{+} \sqcup \mathcal{N}_{x}^+)$ comes a unique minimizer.  
\end{teo}

\section{Structure of  $\Cut \cap \{x=0\}$} \label{sec:x=0}
Describe intersection of the cut locus with the subspace 
  $$M_x = \{q\in M \mid x=0\}.$$

\begin{teo}
There is a stratification
\begin{equation}\label{CutMx}
\Cut \cap M_x = \bigsqcup_{i\in \{+,-\}} \bigg(\mathcal{I}_{xi} \sqcup \mathcal{N}_x^i \sqcup \bigsqcup_{j\in \{+,-\}} \Big(\mathcal{CI}_{xi}^j \sqcup \mathcal{CN}_{xi}^j\Big) \sqcup \mathcal{I}_{zi}^0 \sqcup \mathcal{E}_i \bigg),
\end{equation}
where $\mathcal{I}_{xi} = \bigsqcup_{j\in \{+,-,0\}} \mathcal{I}_{xi}^j$, \ $\mathcal{I}_{xi}^- = \varepsilon^1 (\mathcal{I}_{xi}^+)$, \ $\mathcal{CN}_{xi}^- =\varepsilon^1 (\mathcal{CI}_{xi}^+)$, \ $\mathcal{CN}_{xi}^- = \varepsilon^1(\mathcal{CN}_{xi}^+)$, \ $i\in \{+,-\}$; \ $\mathcal{N}_x^- = \varepsilon^1 (N_x^+)$. Moreover,
\begin{align}
\mathcal{I}_{x+} &= \{q \in M \mid x=0, \ y>0, \ w> G_2 (z,y)\} \cong \R^3, \label{Ix+}\\
\mathcal{I}_{x-} &= \varepsilon^4 (\mathcal{I}_{x+})\cong \R^3 \nonumber\\
\mathcal{N}_x^{\pm} &= \{q \in M \mid x=0, \ \sgn z =\pm 1, \nonumber\\ 
						&\qquad \qquad \qquad -G_3 (z,-y)<w < G_3 (z,y)\} \cong \R^3, \label{Nx+-}
\end{align}
where $G_2, G_3$ are continuous in the set   $\{(z,y) \in \R^2 \mid y>0\}$ functions that satisfy the properties  
\begin{align*}
G_2 (0,y) &= 0,  &G_2(-z,y) = G_2(z,y), \qquad &G_2(\rho^2 z, \rho y) = \rho^3 G_2(z,y), \\ 
&&G_3(-z,y) = G_3 (z,y), \qquad &G_3 (\rho^2 z, \rho y) = \rho^3 G_3 (z,y),  \qquad \rho>0. 
\end{align*}
\end{teo}

\begin{proofn}
By virtue of the equalities $\varepsilon^1 (\Cut) = \Cut$, $\varepsilon^1 (M_{0+})=M_{0-},$ we get, with account of Theorem~\ref{thm:x=0+}, that
$$\Cut \cap M_{0-} = \varepsilon^1 (\Cut \cap M_{0+}) = \bigsqcup_{i\in\{+,-\}} \Big(\mathcal{I}_{xi}^- \sqcup \mathcal{CI}_{xi}^- \sqcup \mathcal{CN}_{xi}^-\Big).$$
Whence, with account of Theorem~\ref{thm:x=z=0}, we obtain stratification~(\ref{CutMx}). Representation~(\ref{Ix+}) is obtained from~(\ref{I+x+}), Table~\ref{tab:dec}, and the equalities  
$$
\mathcal{I}_{x+}^- = \varepsilon^1 (\mathcal{I}_{x+}^+) = \left\{q \in M \mid x=0, z<0, y>0, w>W_{\conj}^{21} (y/\sqrt{|z|}) \sqrt{|z|^3}\right\}
$$
for the function
$$
G_2 (z,y) = \left\{\begin{array}{l}
			W_{\conj}^{21} (y/\sqrt{|z|}) |z|^{3/2}, \  z\neq0, \\
			0, \ z=0.
			\end{array}\right.
$$
Continuity of the function
  $G_2$ on the set   $\{(z,y) \mid z \neq 0, y>0\}$ follows from continuity of the function   $W_{\conj}^{21}$ on the ray $(0,+\infty)$, see Lemma~\ref{lemcon21}. In order to prove continuity of the function   $G_2$ on the ray $\{(z,y) \mid z = 0, y>0\}$, take any sequence    $(z_n,y_n)$, $z_n \to +0$, $y_n \to \bar{y} > 0$. Then, with account of Lemma~\ref{lemcon21}, we get $G_2 (z_n,y_n) = W_{\conj}^{21} (y_n/\sqrt{z_n}) z_n^{3/2} \to 0 = G_2 (0,\bar{y})$. Thus the function $G_2$ is continuous for $y>0$.

Representation~(\ref{Nx+-}) follows from Theorem~\ref{thm:x=0+} for the function $G_3 (z,y) = W_{\conj}^{22} (y/\sqrt{|z|}) |z|^{3/2}.$

Representations~(\ref{Ix+}), (\ref{Nx+-}) imply that the strata   $\mathcal{I}_{x\pm}, \mathcal{N}_x^{\pm}$ are homeomorphic to $\R^3$. The theorem is proved.
\end{proofn}

\section{Global stratification of the cut locus}\label{sec:cut_glob}
In this section we unite results of Sections~\ref{sec:x=z=0}, \ref{sec:z=0}, and \ref{sec:x=0}, and provide a global description of the cut locus. 

\begin{figure}[ht]
\includegraphics[width=0.325\linewidth]{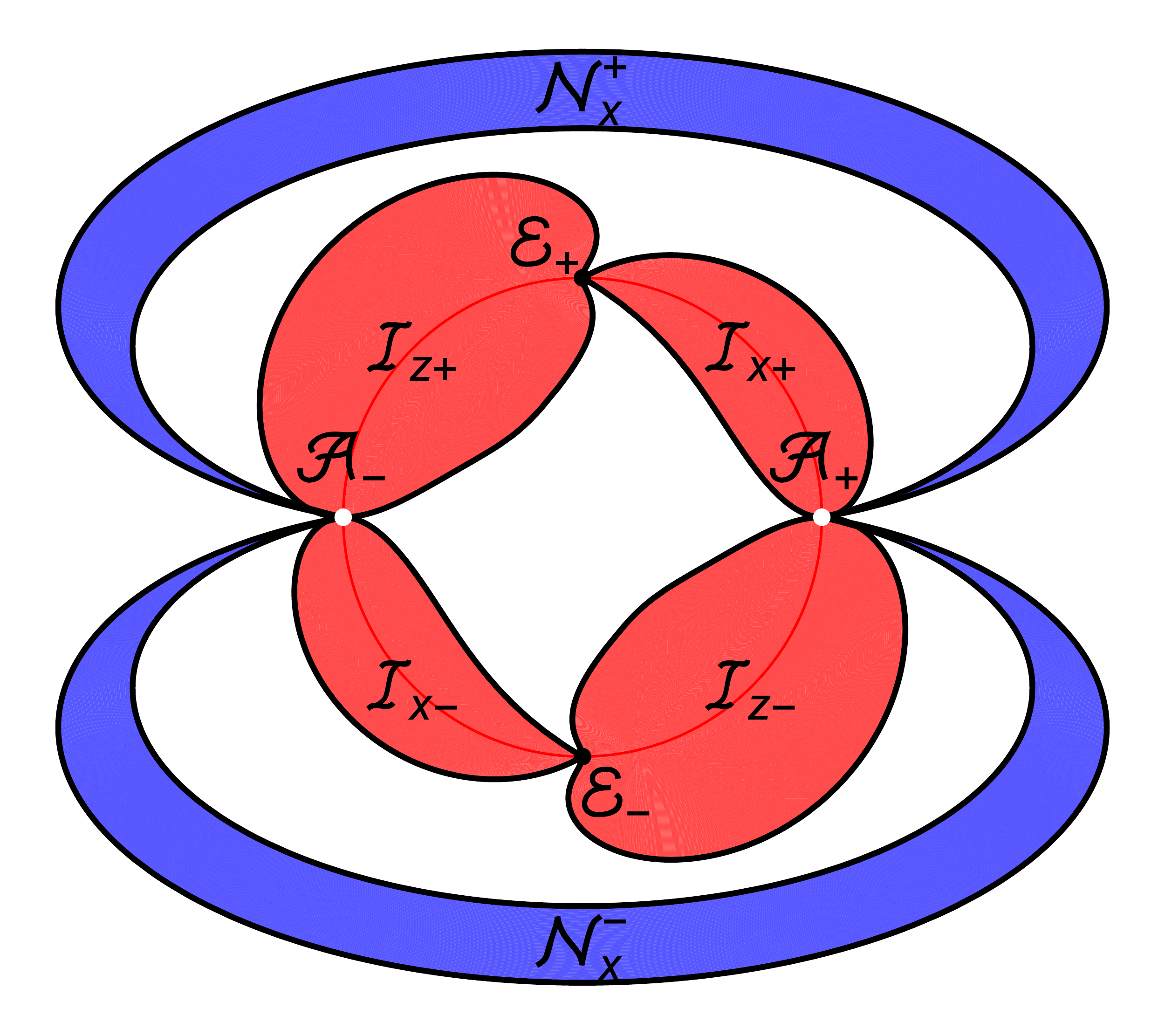} \includegraphics[width=0.325\linewidth]{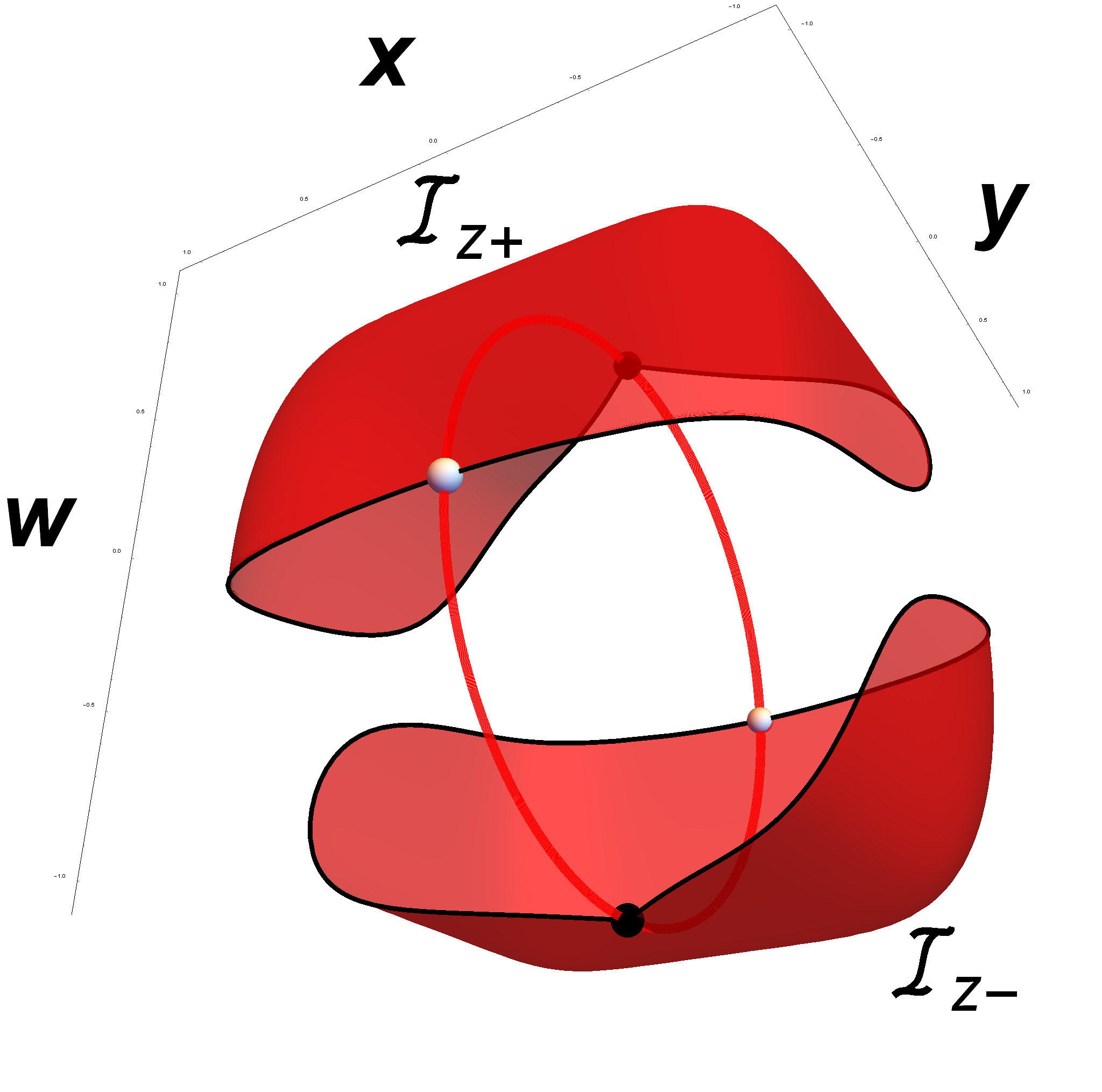} \includegraphics[width=0.325	\linewidth]{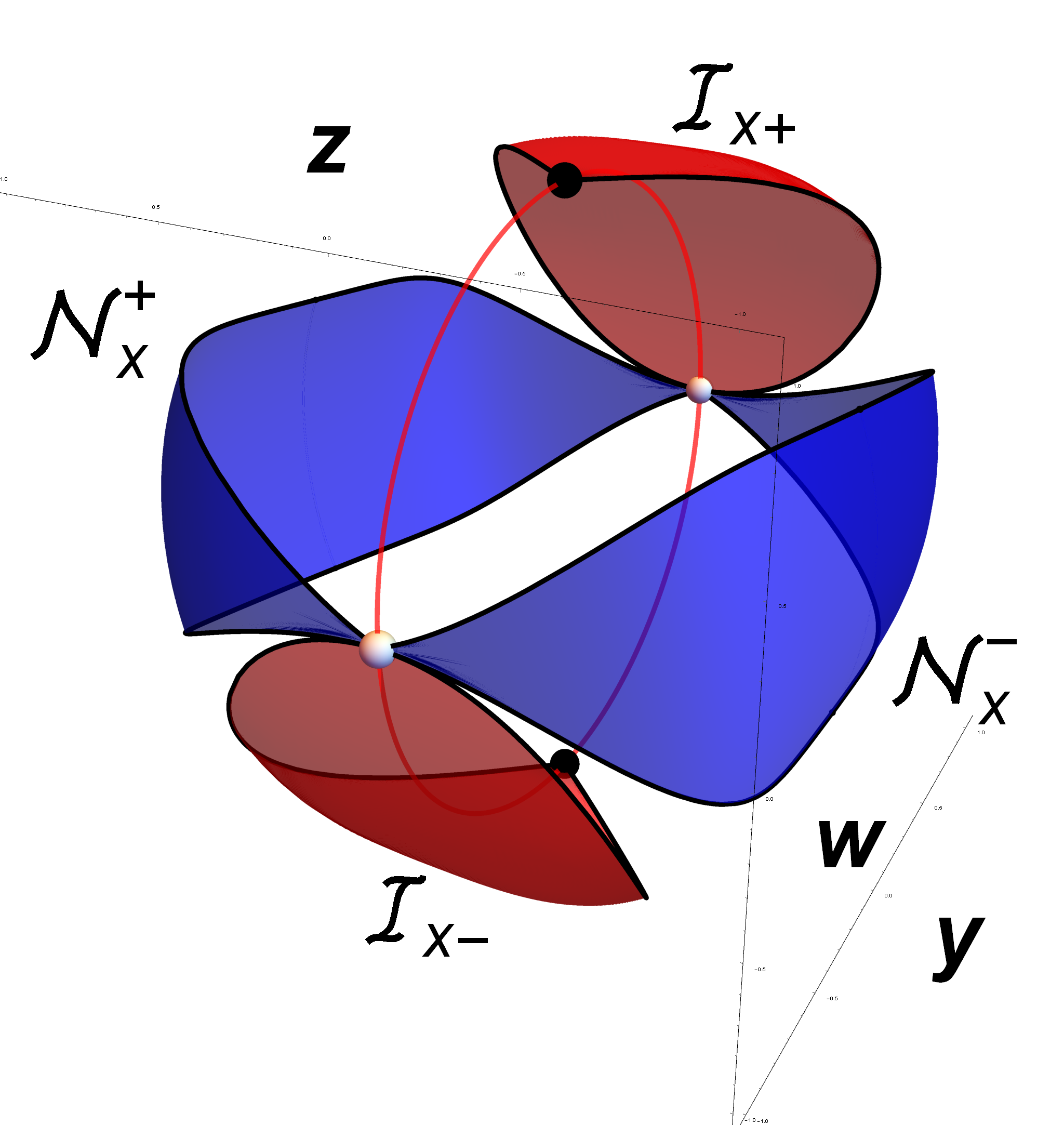}
\caption{Stratification of cut locus: global structure, intersections with the subspaces $M_z$ and $M_x$}\label{fig:dec}
\end{figure}

Figure~\ref{fig:dec} shows stratifications of the cut locus and of its intersections with the coordinate subspaces $M_x$, $M_z$. On the left we show the contiguity topology of strata of the cut locus in the quotient by dilations   $X_0$. In the center we show the set  $\Cut \cap M_z$ after factorization by dilations  $X_0$; the quotient  $M_z / e^{\R X_0}$  is represented by the topological sphere  $\{q\in M \mid x^6 + y^6 + w^2 = 1\}$. Similarly, on the right we show the quotient   $(\Cut \cap M_x) / e^{\R X_0}$ on the topological sphere  $\{q\in M \mid y^6 + |z|^3 + w^2 = 1\}$. 

\begin{teo}\label{thm:global}
The cut locus stratifies as follows:
\begin{align}
\Cut = \bigsqcup_{i \in \{+,-\}} \bigg(\mathcal{I}_{zi} \sqcup \mathcal{I}_{xi} \sqcup \mathcal{N}_x^{i} \sqcup \Big(\bigsqcup_{j\in\{+,-\}} \mathcal{CI}_{zi}^j \sqcup \mathcal{CI}_{xi}^j \sqcup \mathcal{CN}_{xi}^j\Big) \sqcup \mathcal{E}_i\bigg).
\end{align}
Three-dimensional strata
  $\mathcal{I}_{zi}$, $\mathcal{I}_{xi}$, $\mathcal{N}_x^i$, \ $i \in \{+,-\}$, are Maxwell strata, to each point of the strata come two minimizers. Two-dimensional strata   $\mathcal{CI}_{zi}^j$, $\mathcal{CI}_{xi}^j$, $\mathcal{CN}_{xi}^j$, \ $i,j\in \{+,-\}$, consist of conjugate points that are limit points for Maxwell points; to each point of the strata comes a unique minimizer. One-dimensional strata   $\mathcal{E}_i, i \in \{+,-\},$ consist of Maxwell points that are conjugate points; to each point of the strata comes a one-parameter family of minimizers.  
\end{teo}

The cut locus is not closed since it contains points arbitrarily close to the initial point $q_0$, but does not contain the point itself (this is a general fact of sub-Riemannian geometry). Closure of the cut locus  in the sub-Riemannian problem on the Engel group admits the following simple description.

\begin{teo}
$\cl (\Cut) = \Cut  \sqcup \A_+ \sqcup \A_- \sqcup \{q_0\} $.
\end{teo}

Denote by 
$\Conj$ the caustic, i.e., the set of conjugate points along all geodesics starting from the point $q_0$~\cite{engel_conj}; and by $\Max$ the Maxwell set~\cite{engel_cut}. From Th.~\ref{thm:global} we get the following description of the sets $\Cut \cap \Conj$ and $\Cut \cap \Max$.

\begin{teo}
There are stratifications
\begin{align*}
\Cut \cap \Conj &= \bigsqcup_{i \in \{+,-\}, \ j\in\{+,-\}} \Big(\mathcal{CI}_{zi}^j \sqcup \mathcal{CI}_{xi}^j \sqcup \mathcal{CN}_{xi}^j\Big) \sqcup \mathcal{E}_+\sqcup \mathcal{E}_-, \\
\Cut \cap \Max &= \bigsqcup_{i \in \{+,-\}} \bigg(\mathcal{I}_{zi} \sqcup \mathcal{I}_{xi} \sqcup \mathcal{N}_x^{i} \sqcup \mathcal{E}_i\bigg).
\end{align*}
\end{teo}

In other words, $\Cut \cap \Conj$ consists of all two-dimensional and one-dimen-sional strata, while  $\Cut \cap \Max$ consists of all three-dimensional and one-dimensional strata of the cut locus.

\section{Conclusion}

This paper continues the first in the world literature detailed study~\cite{engel, engel_conj, engel_cut} of a sub-Riemannian structure of step more than two. 
A complete description of the cut locus and the structure of optimal synthesis for the nilpotent sub-Riemannian problem on the Engel group is obtained.

Via nilpotent approximation, the results obtained are important for the study and applications of general sub-Riemannian structures of rank 2 in 4-dimensional space. Theoretically, our results open the way to investigation of basic local properties of sub-Riemannian distance for such structures near the initial point. From the applied point of view, these results lead to algorithms and software for solving a path-planning problem in mobile robotics. 

Our research is based upon detailed study of notions and properties introduced and developed by V.\,I.~Arnold~\cite{arnold_osob}:
Maxwell strata and singularities of  a Lagrange mapping generated by a variational problem (exponential mapping of the sub-Riemannian problem). It provides an example of thorough theoretical research related to important applications (Euler's elasticae and mobile robots).  

\section*{APPENDIX~A.}
\setcounter{equation}{0}
\renewcommand{\theequation}{{\rm A}.\arabic{equation}}
In this section we present proofs of the diffeomorphic property for the restriction of the exponential mapping to Maxwell sets for the case 
  $x = z = 0$ considered in Sec.~\ref{sec:x=z=0}.

Denote the functions
\begin{equation}\label{iota12}
\iota_1 (k) = 2E(k)-K(k), \quad  \iota_2 (k) = K(k)-E(k), \qquad k \in (0,1).
\end{equation}

\begin{rem}\label{lemio1}
There hold the inequalities
$\iota_1 (k) > 0$ for $k \in (0,k_0)$, and
$\iota_2 (k) > 0$ for $k \in (0,1)$.
\end{rem}

Introduce an equivalence relation useful for our proofs.
 
\begin{dfn}
Let $X$ be a topological space. Let
 $f_1, f_2 \colon X \to \R, \{\nu_n\} \subset X$. 

We write $f_1 \newtilde f_2$ if $\ds \lim_{n\rightarrow\infty} \frac{f_1 (\nu_n)}{f_2(\nu_n)} \in \R \backslash \{0\}$. 
\end{dfn}

Below in the proofs of the diffeomorphic property of mappings we apply the following Hadamard global diffeomorphism theorem.

\begin{teo}[\cite{Hadamard}]  \label{Hadamard}
Let $F \colon X \to Y$ be a smooth mapping between manifolds of equal dimension. Let the following conditions hold:  
\begin{enumerate}\itemsep=-2pt
\item [$(1)$] $X$ is connected,
\item [$(2)$] $Y$ is connected and simply connected,
\item [$(3)$] $F$ is nondegenerate,
\item [$(4)$] $F$ is proper \big($F^{-1} (K) \subset X$ is compact for a compact $K \subset Y$\big).
\end{enumerate}
Then $F$ is a diffeomorphism.
\end{teo}

\begin{dfn}
A sequence $\{x_n\}$ in a topological space $X$ tends to the boundary of $X$ if there is no compact in   $X$ that contains this sequence.  

Notation: $x_n \to \partial X$.
\end{dfn}

It is easy to see that a continuous mapping 
  $F: X \to Y$ between topological spaces is proper iff for any sequence   $\{x_n\} \subset X$ there holds the implication:    $x_n \to \partial X \Rightarrow F(x_n) \to \partial Y$.

\begin{proofn}[of Lemma~\ref{lemMAX20}]
We apply Th.~\ref{Hadamard}. It follows from definitions of the sets   $\MAX_{++}^{20}$, $\mathcal{I}_{x+}^0$ that $\MAX_{++}^{20}$ is connected, $\mathcal{I}_{x+}^0$ is connected and simply connected, i.e., conditions   (1) and (2) hold for the restriction of   $\Exp$ under consideration. It was shown in work~\cite{engel_conj} that in the case   $u_1 = \pi$ the exponential mapping is nondegenerate for  $\sin u_2 \neq 0$, thus condition~(3)  of Th.~\ref{Hadamard} holds as well.

Let $\nu = (k, u_1,u_2, \sigma) \in \MAX_{++}^{20}$, then
$$\Exp(\nu) = \bigg(0, \frac{4 \iota_1 (k)}{\sigma}, \ 0, \ \frac{8}{3 \sigma^3}\Big(k^2 \iota_1(k) + \iota_2(k)\Big)\bigg).
$$
It follows from Remark~\ref{lemio1} that $\Exp(\nu) \in \mathcal{I}_{x+}^0$, thus $\Exp(\MAX_{++}^{20}) \subseteq \mathcal{I}_{x+}^0$. 

In order to prove condition
  (4), consider any sequence   $\big\{\nu_n = (k_n,\pi,\pi/2,\sigma_n)\big\}$, $n = 1,2,3,\dots$, tending to the boundary of the set   $\MAX_{++}^{20}$ for $n \to \infty$. Denote $\Exp(\nu_n) = (0,y_n,0,w_n)$. We show that $\big\{\Exp(\nu_n)\big\}$ tends to the boundary of   $\mathcal{I}_{x+}^0$. Let us study the possible cases as   $n \to \infty$.  After passing to a subsequence, only the following cases are possible: 
\begin{enumerate}
\item $k_n \to 0$. Then $y_n \newtilde 1/\sigma_n$, whence $w_n \newtilde k_n^2 y_n^3$. Thus $w_n \to 0$ or $y_n \to \infty$, may be, on a subsequence, in the both cases   $\Exp (\nu_n) \to \partial I_{x+}^0$. Below, for brevity, in similar arguments we omit such a phrase about subsequence. 
\item $k_n \to k_0$. Then $y_n \to 0$ or $w_n \to \infty$.
\item $\sigma_n \to \infty, k_n \to \bar{k} \in (0,k_0)$. Then $y_n \to 0, w_n \to 0$.
\item $\sigma_n \to 0, k \to \bar{k} \in (0,k_0)$. Then $y_n \to \infty, w_n \to \infty$.
\end{enumerate}
Thus $\Exp \colon \MAX_{++}^{20} \to \mathcal{I}_{x+}^0$ is a proper mapping, so a diffeomorphism by Th.~\ref{Hadamard}.
\end{proofn}

\begin{proofn}[of Lemma~\ref{lemMAX10}]
It follows from definition of the sets
  $\MAX_{++}^{10}$, $\mathcal{I}_{z+}^0$ that the set   $\MAX_{++}^{10}$ is connected, while   $\mathcal{I}_{z+}^0$ is connected and simply connected, i.e., conditions   (1) and (2) of Th.~\ref{Hadamard} for the restriction of   $\Exp$ under consideration hold.

Let $\nu = (k, u_1,u_2, \sigma) \in \MAX_{++}^{10}$, then 
$$ 
\Exp(\nu) = \bigg(0, \ \frac{2 (2 E_1 -  F_1)}{\sigma}, \ 0, \ \frac{4 E_1 c_1 - d_1^3 s_1}{3 \sigma^3 c_1}\bigg), 
$$
where we used the equality
  $\ds F_1 = 2 E_1 - \frac{s_1 d_1}{c_1}$ equivalent to   $f_z(u_1,k) = 0$. Since $u_1 = u_{1z}(k) \in (\pi/2,\pi)$, then $s_1>0, c_1<0, \  E_1>0, \ 2 E_1 - F_1 < 0,$ whence $\Exp(\nu) \in \mathcal{I}_{z+}^0$, thus $\Exp(\MAX_{++}^{10}) \subseteq \mathcal{I}_{z+}^0$. It was shown in work~\cite{engel_conj} that in the case   $u_1 = u_{1z}(k)$ the exponential mapping is nondegenerate for   $\cos u_2 \neq 0$, thus condition (3) of Th.~\ref{Hadamard} holds.

For the proof of condition
 (4) consider any sequence   $\big\{\nu_n = (k_n,u_{1z}(k_n),0,\sigma_n)\big\}$, $n = 1,2,3,\dots$, tending to the boundary of the set   $\MAX_{++}^{10}$ for $n \to \infty$. Denote $\Exp(\nu_n) = (0,y_n,0,w_n)$. We show that $\big\{\Exp(\nu_n)\big\}$ tends to the boundary of $\mathcal{I}_{z+}^0$. Consider the possible cases as   $n \to \infty$:
\begin{enumerate}
\item $k_n \to k_0$. Then $u_{1z}(k_n) \to \pi$. Whence $y_n \to 0$ or $\sigma_n \to 0$ and $w_n \to \infty$.
\item $k_n \to 1$. Then $u_{1z}(k_n) \to \pi/2$. Whence $y_n \to -\infty$ or $\sigma_n \to \infty$ and $w_n \to 0$.
\item $\sigma_n \to \infty, k_n \to \bar{k} \in (k_0,1)$. Then $y_n \to 0, w_n \to 0$.
\item $\sigma_n \to 0, k_n \to \bar{k} \in (k_0,1)$. Then $y_n \to -\infty, w \to +\infty$.
\end{enumerate}  
So $\Exp \colon \MAX_{++}^{10} \to \mathcal{I}_{z+}^0$ is a proper mapping (condition (4) holds), thus a diffeomorphism by Th.~\ref{Hadamard}.
\end{proofn}

\section*{APPENDIX~B.}
\setcounter{equation}{0}
\renewcommand{\theequation}{{\rm B}.\arabic{equation}}
In this section we prove some lemmas from Sec.~\ref{sec:z=0+}.

\begin{proofn}[of Lemma~\ref{diffzc}]
We differentiate the equality
  $f_z\Big(F\big(u_{1z}(k)\big),k\Big) = 0$ by the variable $k$ and get the following expression for the derivative:   
\begin{align}
u_{1z}'(k) = \frac{ c_1 \Big(1 - \frac{E_1 c_1}{s_1 d_1}\Big)}{k(1 - k^2) s_1}. \label{du1z}
\end{align}
Further, we use this expression to compute derivative of each coordinate of the restricted exponential mapping:
\begin{align}
\frac{\operatorname{d} Y_1^1 \big(k, u_{1z}(k), \pi/2\big)}{\operatorname{d} k} &= \frac{E_1 d_1}{2 k^2 (1 - k^2)^{3/2} s_1} > 0, \label{dY}\\
\frac{\operatorname{d} W_1^1 \big(k, u_{1z}(k), \pi/2\big)}{\operatorname{d} k} &= -\frac{(E_1 c_1 - s_1 d_1^3) (E_1 c_1^3 - s_1 d_1^3)}{16 k^4 (1 - k^2)^{5/2} s_1^6 c_1} > 0,	\label{dW}
\end{align}
since $u_1 = u_{1z}(k) \in (\pi/2,\pi)$ and $E_1>0, s_1>0, c_1<0, d_1>0$. So the mapping   $\Exp\colon \CMAX_{+}^{1+} \to \mathcal{CI}_{z+}^{+}$ is nondegenerate, thus condition (3) of Th.~\ref{Hadamard} holds. Conditions (1) and (2) are also obviously satisfied.   

For the proof of condition (4) consider a sequence 
  $k_n, \ n=1,2,3,\dots,$ tending to the boundary of the set   $\CMAX_{+}^{1+}$ and show that one of the functions   $Y_1^1 \big(k_n, u_{1z}(k_n), \pi/2\big), W_1^1 \big(k_n, u_{1z}(k_n), \pi/2\big)$ tends to infinity:
\begin{enumerate}
\item $k_n \to k_0$, then $u_{1z}(k_n) \to \pi$. We have 
\begin{align}
Y_1^1 \big(k_n, u_{1z}(k_n), \pi/2\big) &\to \frac{1 - 2 k_0^2}{2 k_0 \sqrt{1 - k_0^2}}, \label{Y1-} \\
W_1^1 \big(k_n, u_{1z}(k_n), \pi/2\big) &\newtilde - \frac{E (k_0^2)}{24 k_0^3 (1 - k_0^2)^{3/2} s_1^3} \to -\infty. \label{W1-}
\end{align}
\item $k_n \to 1$, then $u_{1z}(k_n) \to \pi/2$. 
\begin{align}
Y_1^1 \big(k_n, u_{1z}(k_n), \pi/2\big) &= \frac{k_n^2 c_1^2-1 + k_n^2 }{2 k_n \sqrt{1-k_n^2} c_1} \newtilde \frac{k_n^2 c_1}{\sqrt{1 - k_n^2}}-\frac{\sqrt{1-k_n^2}}{c_1}. \label{Y11}
\end{align}
It follows from $f_{z} (F_1,k_n) = 0$ that $\ds c_1^2 = \frac{s_1^2 (1 - k_n^2 s_1^2)}{(F_1-2 E_1)^2},$ compute further
\begin{align*}
\frac{1-k_n^2 s_1^2}{1-k_n^2} &= 1 + \frac{k_n^2 c_1^2}{1-k_n^2} = 1 + \frac{k_n^2 s_1^2 (1-k_n^2 s_1^2 )}{(1-k_n^2)(F_1-2E_1)^2}.
\end{align*}
Since $F_1 - 2 E_1 \to \infty$, we get in the limit  
\begin{align}
\frac{1-k_n^2 s_1^2}{1-k_n^2} = \frac{1}{1-\frac{k_n^2 s_1^2}{(F_1-2 E_1)^2}} \to 1, \label{limd}
\end{align}
consequently,
\begin{align}
\frac{c_1^2}{1-k_n^2} &= \frac{s_1^2 (1-k_n^2 s_1^2)}{(1-k_n^2)(F_1-2E_1)^2} \newtilde \frac{1}{(F_1-2E_1)^2} \to 0. \label{limcos}
\end{align}
From~(\ref{Y11}) we get 
\begin{align}
Y_1^1 \big(k_n, u_{1z}(k_n), \pi/2\big) \newtilde \frac{\sqrt{1-k_n^2}}{c_1} \to \infty. \label{Y1+}
\end{align} 

Consider the second coordinate:
\begin{align}
W_1^1 \big(k_n, u_{1z}(k_n), \pi/2\big) &= -\frac{d_1^3 E_1}{48 k_n^3 (1 - k_n^2)^{3/2} s_1^3} -\frac{\sqrt{1 - k_n^2} (-1 + 4 k_n^2)}{48 k_n^3  c_1}  \nonumber\\
&\quad +\frac{c_1 \big( c_1^2 + (1-k_n^2) (1 +k_n^2+4k_n^4) s_1^2 \big)}{48 k_n^3 (1 - k_n^2)^{3/2} s_1^2}.\label{W11}
\end{align}
Further consider~(\ref{W11}) term by term, using~(\ref{limd}) and (\ref{limcos}): 
\begin{align*}
-\frac{d_1^3 E_1}{48 k_n^3 (1 - k_n^2)^{3/2} s_1^3 } &\to - \frac{1}{48},\qquad &-\frac{\sqrt{1 - k_n^2} (-1 + 4 k_n^2)}{48 k_n^3  c_1} &\to \infty,\\
\frac{c_1^3}{48 k_n^3 (1 - k_n^2)^{3/2} s_1^2 } &\to 0, &\frac{c_1(1 +k_n^2+4k_n^4)}{48 k_n^3 \sqrt{1 - k_n^2}} &\to 0, 
\end{align*}
whence it follows that
\begin{align}
W_1^1 \big(k_n, u_{1z}(k_n), \pi/2\big) \newtilde \frac{\sqrt{1-k_n^2}}{c_1} \to \infty. \label{W1+}
\end{align}
\end{enumerate}
Thus $\Exp\colon \CMAX_{+}^{1+} \to \mathcal{CI}_{z+}^{+}$ is a proper mapping, so a diffeomorphism by Th.~\ref{Hadamard}.
\end{proofn}

\begin{proofn}[of Lemma~\ref{lem:c-z=0+}]
\begin{itemize}
\item[1)] follows from equalities~(\ref{dY})--(\ref{W1-}), (\ref{Y1+}), (\ref{W1+}). 
\item[2)] follows from the inequality   $$ Y_1^1 (k, u_1, \pi/2) - 6 W_1^1 (k, u_1, \pi/2) = \frac{d_1^3 (c_1 E_1 - d_1 s_1)}{8 k^3 (1 - k^2)^{3/2} s_1^3 c_1} > 0$$ under the condition $k\in (k_0,1), \ u_1 \in (\pi/2,\pi)$.
\item[3)] follows from (\ref{Y1+}), (\ref{W1+}).
\end{itemize}
\end{proofn}

Now we prove a lemma we will use in the sequel for localization of two-dimensional parametrically defined sets.

\begin{lem}\label{below}
Consider a parametrically defined set in the plane $(X_1,X_2)\in\R^2$ 
\begin{align*}
\Omega = \Big\{\big(X_1,X_2\big) = \big(f_1(x_1,x_2),f_2(x_1,x_2)\big) \mid x_1 \in (x_1^0,x_1^1), x_2 \in (x_2^0, x_2^1)\Big\},
\end{align*}
where $f_1,f_2 \in C^1\big((x_1^0,x_1^1)\times[x_2^0,x_2^1)\big)$. Consider also a curve 
\begin{align*}
\gamma_0 = \Big\{\big(X_1,X_2\big) = \big(f_1(x_1,x_2^0),f_2(x_1,x_2^0)\big) \mid x_1 \in (x_1^0, x_1^1)\Big\}.
\end{align*}
Let the curve $\gamma_0$ divide the plane $(X_1,X_2)$ into two connected components, and let the following conditions hold for   $x_1 \in (x_1^0,x_1^1)$:
\begin{align}
&\frac{\partial f_1}{\partial x_1}(x_1,x_2)>0, \qquad &x_2\in[x_2^0,x_2^1), \label{partialx1}\\
&\frac{\partial f_1}{\partial x_2}(x_1,x_2)>0, \qquad  &x_2\in(x_2^0,x_2^1), \label{partialx2} \\
&\frac{\partial f_1}{\partial x_2}(x_1,x_2^0)\geq 0, \qquad & \label{partialx20}\\
&\nabla (x_1,x_2) = \bigg(\frac{\partial f_2/\partial x_1}{\partial f_1/\partial x_1}-\frac{\partial f_2/\partial x_2}{\partial f_1/\partial x_2}\bigg)\big(x_1,x_2\big) > 0, &x_2\in[x_2^0,x_2^1), \label{nabla} 
\end{align}
then condition~(\ref{partialx1}) allows us to invert the function  $X_1 = f_1(x_1,x_2^0)$ in the interval $(x_1^0,x_1^1)$: \quad
$x_1 =  \mathbf{h}_0(X_1)$, 
which allows us to define the function 
  $\gamma_0$ as a graph:
$X_2= \mathbf{g}_0(X_1)$,
where $\mathbf{g}_0(X_1) = f_2 (\mathbf{h}_0(X_1),x_2^0)$. 

Then there holds the following condition:  
\begin{align}
f_2(x_1,x_2) < \mathbf{g}_0\big(f_1(x_1,x_2)\big), \qquad x_1 \in (x_1^0,x_1^1), \ x_2 \in (x_2^0, x_2^1), \label{condunder}
\end{align}
i.e., the set
  $\Omega$ lies below the curve $\gamma$ in the plane $(X_1, X_2)$.
\end{lem}
\begin{proofn}
The set $\Omega$ is a union of the curves
$$\gamma_a=\Big\{\big(X_1,X_2\big) = \big(f_1(x_1,x_2^a),f_2(x_1,x_2^a)\big) \mid x_1 \in \big(x_1^0,x_1^1\big)\Big\},$$ 
where $x_2^a=x_2^0 + (x_2^1-x_2^0)a$ for $a \in (0, 1)$. Condition~(\ref{partialx1})  allows us to define as a graph not only the curve   $\gamma_0$, but each curve from the family  $\{\gamma_a \mid a \in (0, 1)\}$ as well:
\begin{align*}
X_2 &= \mathbf{g}_a(X_1), \qquad X_1 \in \big(f_1(x_1^0,x_2^a), f_1(x_1^1,x_2^a)\big),
\end{align*}
where $\mathbf{g}_a(X_1) = f_2 (\mathbf{h}_a(X_1),x_2^a)$, and $x_1 = \mathbf{h}_a(X_1)$  is the inverse function to   $X_1 = f_1(x_1,x_2^a)$.

On the other hand, by fixing 
 $x_1$ the set $\Omega \cup \gamma_0$ becomes a union of curves defined on a half-interval:  
\begin{align*}
\beta_k = \Big\{\big(X_1,X_2\big) = \big(f_1(x_1^k,x_2),f_2(x_1^k,x_2)\big) \mid x_2 \in [x_2^0,x_2^1)\Big\}, 
\end{align*}
where $x_1^k=x_1^0 + (x_1^1-x_1^0)k$ for $k \in (0, 1)$. Moreover, conditions~(\ref{partialx2}), (\ref{partialx20}) allow us to define as a graph each curve of the family   $\{\beta_k \mid k \in (0, 1)\}$:
\begin{align*}
X_2 = \mathbf{v}_{k}(X_1), \qquad X_1 \in \big[f_1(x_1^k,x_2^0), f_1(x_1^k,x_2^1)\big),
\end{align*}
where $\mathbf{v}_k(X_1) = f_2 (x_1^k,\mathbf{r}_k(X_1))$, and $x_2 = \mathbf{r}_k(X_1)$  is the inverse function to   $X_1 = f_1(x_1^k,x_2)$.

Notice that since   $f_1,f_2 \in C^1\big((x_1^0,x_1^1)\times[x_2^0,x_2^1)\big)$, then by inverse function theorem conditions (\ref{partialx1})--(\ref{nabla}) guarantee that   
$\mathbf{g}_a \in C^1 \Big(\big(\bar{f}_1(x_1^0,x_2^a), \bar{f}_1(x_1^1,x_2^a)\big)\Big)$, $\mathbf{v}_k \in C^1\Big(\big[f_1(x_1^k,x_2^0), \bar{f}_1(x_1^k,x_2^1)\big)\Big)$, 
moreover, condition~(\ref{partialx1}) implies existence of a limit   $\lim\limits_{x_1 \to x_1^0} f_1(x_1,x_2^a) =: \bar{f}_1(x_1^0,x_2^a) \in [-\infty,+\infty]$, we define similarly  $\bar{f}_1(x_1^1,x_2^a), \bar{f}_1(x_1^k,x_2^1)$.

Let $\check{X}_1 = f_1 (x_1^k, x_2^a)$, then we get   $\nabla (x_1^k,x_2^a) =  \mathbf{g}_a'(\check{X}_1) - \mathbf{v}_k'(\check{X}_1) > 0$, whence it follows that
\begin{align*}
\forall x_1^k \in (x_1^0,x_1^1) \ x_2^a \in [x_2^0,x_2^1) \ \  \exists \epsilon > 0 \ \ \forall X_1 \in (\check{X_1}, \check{X_1}+\epsilon) \qquad \mathbf{g}_a(X_1) > \mathbf{v}_k(X_1). 
\end{align*}
If $a=0$, then definitions of the functions   $\mathbf{g}_a, \mathbf{v}_k$ imply that
\begin{align*}
\forall x_1^k \in (x_1^0,x_1^1) \ \  \exists \delta > 0 \ \ \forall x_2 \in (x_2^0, x_2^0+\delta) \qquad \mathbf{g}_0\big(f_1 (x_1^k,x_2)\big ) > f_2 (x_1^k,x_2). 
\end{align*}

Suppose that condition~(\ref{condunder}) is violated, i.e.,
\begin{align}
\exists \hat{x}_1 \in (x_1^0,x_1^1) \ \ \exists \hat{x}_2 \in (x_2^0, x_2^1) \qquad \mathbf{g}_0 \big(f_1(\hat{x}_1,\hat{x}_2)\big) = f_2(\hat{x}_1,\hat{x}_2). \label{suggest}
\end{align}

Introduce the notation $\hat{X}_i=f_i(\hat{x}_1,\hat{x}_2), i=1,2$, and $\check{x}_1 = \mathbf{h}_0(\hat{X}_1)$. 
Notice that by definition of
  $\mathbf{g}_0$ we have
$
f_i(\check{x}_1,x_2^0)=\hat{X}_i, i=1,2.
$ 
It follows from (\ref{partialx2}), (\ref{partialx20}) that $f_1(\hat{x}_1, x_2^0)<f_1(\hat{x}_1,\hat{x}_2) = \hat{X}_1 = f_1 (\check{x}_1,x_2^0)$, thus condition (\ref{partialx1}) implies that $\hat{x}_1 < \check{x}_1$.

Below we define a function
  $x_2 = \omega(x_1)$ on the segment $x_1 \in [\hat{x}_1, \check{x}_1]$, which satisfies the condition
\begin{align}
f_1 \big(x_1, \omega(x_1)\big) = \hat{X}_1.  \label{f1const}
\end{align}
In view of condition~(\ref{partialx2}), if a function $w$ is defined, then it is unique. At the endpoints of the segment   $[\hat{x}_1, \check{x}_1]$ the function is defined: $\omega(\hat{x}_1) = \hat{x}_2$, $\omega(\check{x}_1) = x_1^0$,
moreover, we have  
\begin{align}
f_2 \big(\hat{x}_1, \omega(\hat{x}_1)\big) = f_2 \big(\check{x}_1, \omega(\check{x}_1)\big) = \hat{X}_2. \label{f2const}
\end{align}

Notice that inequality~(\ref{partialx1}) implies that for   $x_1 \in (\hat{x}_1,\check{x}_1)$ we have
\begin{align*}
f_1(x_1,x_2^0) < f_1(\check{x}_1,x_2^0) = \hat{X}_1 = f_1(\hat{x}_1, \hat{x}_2) < f_1(x_1,\hat{x}_2),
\end{align*}
which, together with continuity of the function
  $f_1$, implies that the function   $\omega$ is defined at the segment   $x_1 \in [\hat{x}_1,\check{x}_1]$. If $x_1^k \in [\hat{x}_1,\check{x}_1]$, then $\omega (x_1^k) = \mathbf{r}_k (\hat{X}_1)$. In other words,   $x_2 = \omega(x_1^k)$ is a uniquely defined function, inverse to   $\hat{X}_1 = f_1(x_1^k,x_2)$ at the segment $x_1 \in [\hat{x}_1,\check{x}_1]$. By the inverse function theorem,   $\omega \in C^1([\hat{x}_1,\check{x}_1])$.

Denote the functions
$\hat{f}_i(x_1) = f_i \big(x_1, \omega(x_1)\big)$, $i = 1, 2$, and
compute their derivatives:
\begin{align*}
\frac{\operatorname{d} \hat{f}_i}{\operatorname{d} x_1} (x_1) = \frac{\partial f_i}{\partial x_1} \big(x_1, \omega(x_1)\big)+ \frac{\partial f_i}{\partial x_2} \big(x_1, \omega(x_1)\big) \frac{\operatorname{d} w}{\operatorname{d} x_1} (x_1), \qquad i =1,2.
\end{align*}
It follows from (\ref{f1const}) that $\ds \frac{\operatorname{d} \hat{f}_1}{\operatorname{d} x_1} = 0$, whence
\begin{align*}
\frac{\operatorname{d} \hat{f}_2}{\operatorname{d} x_1} (x_1) &= \bigg(\frac{\partial f_2}{\partial x_1} + \frac{\partial f_2}{\partial x_2} \Big(- \frac{\partial f_1/\partial x_1}{\partial f_1/\partial x_2}\Big)\bigg) (x_1, \omega(x_1)) 
= \bigg(\frac{\partial f_1}{\partial x_1} \nabla\bigg) \big(x_1, \omega(x_1)\big) > 0,
\end{align*}
where $x_1 \in [\hat{x}_1, \check{x}_1)$, thus the function   $\hat{f}_2(x_1)$ increases at the segment  $[\hat{x}_1, \check{x}_1)$, whence by continuity of the function we come to a contradiction with condition  (\ref{f2const}).Thus assumption~(\ref{suggest})  is violated, q.e.d.  
\end{proofn}

\begin{lem}\label{subset1}
There holds the inclusion
  $\Exp(\MAX_{+-}^{1+}) \subseteq \mathcal{I}_{z+}^+$.
\end{lem}
\begin{proofn}
Notice that for
  $k\in (k_0,1), u_2 \in (\pi/2,\pi)$ we have
\begin{align}
J_{11}^1(k,u_2)&=\frac{\partial Y_1^1 \big(k, u_{1z}(k), u_2 \big)}{\partial u_2} = \frac{\Delta c_2}{2 k c_1 s_2^2 d_2^3} > 0, \label{dYu2}\\ 
J_{21}^1(k,u_2)&=\frac{\partial Y_1^1 \big(k, u_{1z}(k), u_2 \big)}{\partial k} =  \frac{\Delta (E_1 c_1 d_2^2 - k^2 s_1 d_1 c_2^2)}{2 k^2 (1 - k^2) s_1 c_1 d_1 s_2 d_2^3} > 0, \label{dYk}\\
\nabla^1(k,u_2)&= \frac{\partial W_1^1\big(k,u_{1z}(k),u_2\big)/(\partial k)}{\partial Y_1^1 \big(k,u_{1z}(k),u_2\big)/(\partial k)} - \frac{\partial W_1^1\big(k,u_{1z}(k),u_2\big)/(\partial u_2)}{\partial Y_1^1 \big(k,u_{1z}(k),u_2\big)/(\partial u_2)}  \nonumber \\
&= \frac{\Delta^2 \Big(-2 E_1 s_1 c_1 d_1^5   + s_1^2 d_1^6 + E_1^2 c_1^2 \big((1-k s_1^2)^2 + 2 k (1-k) s_1^2 \big)\Big)}
			{8 k^2  s_1^5  d_1^5 s_2^2 (k^2 s_1 d_1 c_2^2-E_1 c_1 d_2^2)}> 0 \label{cond2}.
\end{align}
Moreover, for all
  $k\in (k_0,1)$ there exists $\lim\limits_{u_2 \to \pi/2} \nabla^1(k,u_2) \in (0,+\infty)$, and there holds the inequality $J_{21}^1(k,\pi/2)>0$. Whence, with account of Lemma~\ref{lem:c-z=0+}, it follows that the mapping   $\restr{\Exp}{\MAX_{+-}^{1+}}$ satisfies conditions of Lemma~\ref{below}, consequently, $\Exp(\MAX_{+-}^{1+}) \subseteq \mathcal{I}_{z+}^+$.
\end{proofn}

\begin{proofn}[of Lemma \ref{diffzm}]
Conditions (1), (2) of Theorem~\ref{Hadamard} are obviously satisfied, and results of~\cite{engel_conj} imply that condition (3) holds as well. It remains to check validity of condition (4).  

Consider arbitrary sequences $(k_n,u_{1z}(k_n),u_{2}^n), \ n=1,2,3,\dots$, in the image of
  $\Exp$ tending to boundary of the set    $\MAX_{+-}^{1+}$ and show that either   $\Exp\big(k_n, u_{1z}(k_n), u_{2}^n\big) \to \mathcal{CI}_{z+}^{+}$, or one of coordinates in the image tends to infinity: 
\begin{enumerate}
\item $u_{2}^n\to \pi/2$, by definition of $\mathcal{CI}_{z+}^{+}$ we have $\Exp(k_n, u_{1z}(k_n), u_{2}^n) \to \mathcal{CI}_{z+}^{+}.$
\item $u_{2}^n\to \pi$, then $\ds Y_1^1 \big(k_n,u_{1z}(k_n),u_{2}^n\big) \newtilde - \frac{1}{2 k_n c_1 s_2} \to \infty.$
\item $k_n \to k_0, u_{2}^n \to \hat{u}_2 \in (\pi/2,\pi)$, then $u_{1z}(k_n) \to \pi$, thus   $s_1 \to 0, c_1 \to -1, d_1 \to 1, \Delta \to 1$. We get 
$\ds W_1^1 \big(k_n,u_{1z}(k_n),u_{2}^n\big) \newtilde - \frac{E(k_0)}{24 k_0^3 s_1^3 s_2^3 d_2^3} \to - \infty.$ 
\item $k_n \to 1, u_{2}^n \to \hat{u}_2 \in (\pi/2,\pi)$, then $u_{1z}(k_n) \to \pi/2$. By virtue of $c_1 \to 0, c_2 \neq 0, s_2 \neq 0$, we get
$\ds Y_1^1 \big(k_n,u_{1z}(k_n),u_{2}^n\big) \newtilde - \frac{c_2}{2 c_1 s_2} \to \infty.$
\end{enumerate}
Consequently, $\Exp\colon \MAX_{+-}^{1+} \to \mathcal{I}_{z+}^+$ is proper, with account of Lemma~\ref{subset1}.  Whence by Th.~\ref{Hadamard} the restriction of the exponential mapping is a diffeomorphism.  
\end{proofn}

\section*{APPENDIX~C.}
\setcounter{equation}{0}
\renewcommand{\theequation}{{\rm C}.\arabic{equation}}
In this Appendix we prove some technical lemmas from Section~\ref{sec:x=0+}.

\begin{proofn}[of Lemma~\ref{lem:c2-x=0}]
We apply Th.~\ref{Hadamard}. Notice first that conditions (1), (2) hold. Further compute derivative of each coordinate by   $k$:
\begin{align}
\frac{\operatorname{d}  Y_1^2(k,0)}{\operatorname{d}  k} &=-\frac{E(k)}{(1 - k^2) \sqrt{2 k^3 \iota_1(k)}} < 0, \label{dY21}\\
\frac{\operatorname{d}  W_1^2(k,0)}{\operatorname{d}  k} &=\frac{\iota_3 (k)}{\big(2 k \iota_1(k)\big)^{5/2}},
\end{align}
where $\ds \iota_3(k)=\frac{3 E^2(k) - (5 - 4 k^2) E(k) K(k) + 2 (1 - k^2) K^2(k)}{1-k^2}>0$ for $k \in (0,k_0)$, since $\iota_3(0)=0$, and $\ds\iota_3'(k) = \frac{\iota_2^2(k) +   k^2 \iota_1(k) \big(2 E(k)+K(k)\big)}{k (1 - k^2)^2} > 0$. Thus for   $k \in (0, k_0)$ we have
\begin{align}
\frac{\operatorname{d}  W_1^2(k,0)}{\operatorname{d}  k} > 0. \label{dW21}
\end{align}
Thus the mapping
  $\Exp: \CMAX_{1+}^{2+} \to \mathcal{CI}_{x+}^{+}$ is nondegenerate, and condition (3) of Th.~\ref{Hadamard} holds.

In order to prove condition (4), consider an arbitrary sequence 
  $k_n, n = 1,2,3, \dots,$ tending to the boundary of  $\CMAX_{1+}^{2+}$ and show that one of coordinates tends to infinity:  
\begin{enumerate}
\item $k_n \to 0$, then 
\begin{align}
Y_1^2(k_n,0) \to \infty, \qquad W_1^2(k_n,0) \to 0. \label{lim21-}
\end{align}
\item $k_n \to k_0$, then
\begin{align}
Y_1^2(k_n,0) \to 0, \qquad W_1^2(k_n,0) \to \infty. \label{lim21+}
\end{align}
\end{enumerate}
Whence  $\Exp: \CMAX_{1+}^{2+} \to \mathcal{CI}_{x+}^{+}$ is a proper mapping, thus a diffeomorphism by Th.~\ref{Hadamard}. 
\end{proofn}

\begin{proofn}[of Lemma~\ref{lemcon21}]
Follows from Lemma~\ref{lem:c2-x=0} and expressions~(\ref{dY21})--(\ref{lim21+}).
\end{proofn}

\begin{lem} \label{subset21}
There holds the inclusion
  $\Exp (\MAX_{1++}^{2+}) \subset \mathcal{I}_{x+}^+$.
\end{lem}
\begin{proofn}
Notice the inequalities
\begin{align}
J_{11}^{21}=\frac{\partial Y_1^2(k,u_2)}{\partial u_2} &= \sqrt{\frac{\iota_1(k)}{2 k c_2^3}} s_2 > 0, \label{J1121}\\
J_{12}^{21}=\frac{\partial W_1^2(k,u_2)}{\partial u_2} &= \frac{\big(\iota_2(k) +k^2 s_2^2 \iota_1(k)\big) s_2}{2 c_2 \big(2 k \iota_1(k) c_2\big)^{3/2}} > 0. \label{J1221}
\end{align}
Thus, if 
  $u_2$ grows, then the both coordinates   $Y^2$, $W^2$ grow as well. Thus for a fixed   $k=\hat{k}$ the curve $\big(Y^2(\hat{k},u_2),W^2(\hat{k},u_2)\big), u_2 \in (0,\pi/2)$,  lies above the curve  $W^2=W_{\conj}^{21}(Y^2)$, $Y_2 > Y^2(\hat{k},0)$, see Lemma~\ref{lemcon21}. This proves this lemma by definition~(\ref{Ix++}).
\end{proofn}

\begin{proofn}[of Lemma~\ref{diffxm1}]
Apply Th.~\ref{Hadamard}. Conditions (1), (2) follow from definitions of the sets   $\MAX_{1++}^{2+}$, $\mathcal{I}_{x+}^+$. Condition (3) follows from results of~\cite{engel_conj}.  With account of Lemma~\ref{subset21}, in order to prove condition (4), it suffices to consider arbitrary sequences $(k_n, u_{2}^n)$ in the image of  $\Exp$ tending to boundary of the set   $\MAX_{1++}^{2+}$, and to show that the sequence   $\big(Y_1^2(k_n,u_{2}^n),W_1^2(k_n,u_{2}^n)\big)$ tends to boundary of the set   $\mathcal{I}_{x+}^+$:
\begin{enumerate}
\item $u_{2}^n\to 0,$ by definition of $\mathcal{CI}_{x+}^{+}$ we have $$\big(Y_1^2(k_n,u_{2}^n),W_1^2(k_n,u_{2}^n)\big) \to \mathcal{CI}_{x+}^{+}.$$
\item $u_{2}^n \to \pi/2, k_n \to \hat{k} \in [0,k_0)$, or $k_n \to 0$, then $$Y_1^2(k_n,u_{2}^n) \newtilde \frac{1}{\sqrt{k_n c_2}} \to \infty.$$
\item $u_{2}^n \to \pi/2, k_n \to k_0$, then $\iota_1(k_n) \to 0$, thus $$W_1^2(k_n,u_{2}^n) \newtilde \frac{E(k_0)}{(\iota_1(k) c_2)^{3/2}} \to \infty.$$
\item $k_n \to k_0, u_{2}^n \to \hat{u}_2 \in (0,\pi/2)$, then $\iota_1(k_n) \to 0$, thus $$Y_1^2(k_n,u_{2}^n) \newtilde \sqrt{\iota_1(k_n)} \to 0.$$
\end{enumerate} 
So $\Exp \colon \MAX_{1++}^{2+} \to \mathcal{I}_{x+}^+$ is a proper mapping, thus a diffeomorphism by Th.~\ref{Hadamard}.
\end{proofn}

Introduce the functions
\begin{align}
\iota_4(k) &= (2 - k^2) K(k) - 2 E(k), \label{iota4}\\ 
\iota_5(k) &= (2 - k^2) E(k) K(k) + (1 - k^2) K^2(k) -3 E^2(k), \nonumber\\
\iota_6(k) &= E(k) + (-1 + k^2) K(k) > 0. \nonumber
\end{align}

\begin{lem}
Let $k \in (0,1)$. There hold the inequalities   $\iota_i(k)>0, \ i=4,5,6$.
\end{lem}
\begin{proofn}
Notice that $\iota_i(0) = 0, i=4,5,6$. Further, $$\bigg(\frac{\iota_4(k)}{\sqrt{1 - k^2}}\bigg)' = \frac{k \iota_2(k)}{(1 - k^2)^{3/2}} > 0$$ for $k \in (0,1)$;
since $\big((2 - k^2) E(k) -  2 (1 - k^2) K(k)\big)'=3 k \iota_2(k) > 0$, then 
$$\ds \bigg(\frac{\iota_5(k)}{\sqrt{1 - k^2}}\bigg)' = \frac{\iota_4(k) \big((2 - k^2) E(k) -  2 (1 - k^2) K(k)\big)}{k (1 - k^2)^{3/2}}>0;
$$ 
finally,  $\iota_6'(k)=k K(k)>0.$
Whence $\iota_i(k) > 0, i=4,5,6,$ for $k \in (0,1)$.
\end{proofn}

\begin{proofn}[of Lemma~\ref{diffxc2+}]
Apply Th.~\ref{Hadamard}. Notice that conditions (1), (2) hold. Compute further derivative of each coordinate in   $k$:
\begin{align}
\frac{\partial Y_2^2(k,0)}{\partial k} &= -\frac{k \iota_2(k)}{2 (1 - k^2)^{5/4} \sqrt{\iota_4(k)}} <0, \label{dY22+} \\
\frac{\partial W_2^2(k,0)}{\partial k} &= \frac{k^3 \iota_5(k)}{4 (1 - k^2)^{7/4} \iota_4^{5/2}(k)} >0. \label{dW22+}
\end{align}
By virtue of sign-definiteness of the derivatives, the mapping is nondegenerate, thus condition (3) of Th.~\ref{Hadamard} holds.

In order to prove condition (4), consider arbitrary sequences 
  $k_n, n = 1,2,3,\dots,$ tending to boundary of the set   $\CMAX_{2++}^{2+}$:
\begin{enumerate}
\item $k_n \to 0$, then 
\begin{align}
Y_2^2(k_n,0) \to 0, \qquad W_2^2(k_n,0) \to 1/\sqrt{\pi}. \label{lim22+-}
\end{align}
\item $k_n \to 1$, then
\begin{align}
Y_2^2(k_n,0) \to -\infty, \qquad W_2^2(k_n,0) \to \infty. \label{lim22++}
\end{align}
\end{enumerate}
Whence $\Exp: \CMAX_{2++}^{2+} \to \mathcal{CN}_{x++}^{+}$ is a proper mapping, thus a diffeomorphism by Th.~\ref{Hadamard}.
\end{proofn}

\begin{proofn}[of Lemma~\ref{lemcon22+}]
Follows from Lemma~\ref{diffxc2+} and expressions~(\ref{dY22+})--(\ref{lim22++}).
\end{proofn}

\begin{proofn}[of Lemma~\ref{diffxc2-}]
Apply Th.~\ref{Hadamard}. Conditions (1), (2) hold. Compute derivative of each coordinate in   $k$:
\begin{align}
\frac{\partial Y_2^2(k,\pi/2)}{\partial k} &= -\frac{k \iota_6 (k)}{2 (1 - k^2) \sqrt{\iota_4(k)}} < 0, \label{dY22-}  \\
\frac{\partial W_2^2(k,\pi/2)}{\partial k} &= \frac{k^3 \iota_5(k)}{4 (1 -  k^2) \iota_4^{5/2}(k)} > 0,
\end{align}
thus the mapping $\restr{\Exp}{\CMAX_{2+-}^{2+}}$ is nondegenerate.

Consider arbitrary sequences
  $k_n, n=1,2,3,\dots,$ tending to boundary of the set   $\CMAX_{2+-}^{2+}$:
\begin{enumerate}
\item $k \to 0$, then 
\begin{align}
Y_2^2(k_n,\pi/2) \to 0, \qquad W_2^2(k_n,\pi/2) \to -1/\sqrt{\pi}. \label{lim22--}
\end{align}
\item $k \to 1$, then 
\begin{align}
Y_2^2(k_n,\pi/2) \to -\infty, \qquad W_2^2(k_n,\pi/2) \to 0. \label{lim22-+}
\end{align}
\end{enumerate}
Whence $\Exp: \CMAX_{2+-}^{2+} \to \mathcal{CN}_{x+-}^{+}$ is proper, thus a diffeomorphism by Th.~\ref{Hadamard}.
\end{proofn}

\begin{proofn}[of Lemma~\ref{lemcon22-}]
Follows from Lemma~\ref{diffxc2-} and expressions~(\ref{dY22-})--(\ref{lim22-+}).
\end{proofn}

\begin{proofn}[of Lemma~\ref{conjcross}]
To prove the lemma, consider the following fixed points of the symmetry
  $\varepsilon^2$ for $\lambda \in C_3$: \quad
$$\FIX_{3+}^{2+} = \Big\{ (\sigma, p, \tau) \in C_{3+}^+ \times \R_+ \mid \tau = 0\Big\} \subset \FIX^{2}.$$
The exponential mapping transforms the set
  $\FIX_3^{2+}$ into the set  
\begin{align*}
\Fix_{3+}^{2+} &= \Big\{ (Y^2,W^2) = \big(Y_3^2(p), W_3^2(p)\big) \mid p\in (0,\infty) \Big\},\\
\big(Y_3^2(p),W_3^2(p)\big) &= \bigg(\frac{2 \sinh p - p \cosh p}{\sqrt{p \cosh p  - \sinh p}}, \frac{9 \sinh p - 12 p \cosh p + \sinh (3 p)}{24 (p \cosh p - \sinh p)^{\frac{3}{2}}} \bigg).
\end{align*}

Since $\big(2 \sinh p - p \cosh p\big)' = \cosh p - p \sinh p$, then
\begin{align}
\frac{\operatorname{d} Y_3^2(p)}{\operatorname{d} p} < 0,
\end{align}
moreover,
\begin{align}
\lim_{p\to 0}Y_3^2(p) = +\infty, \qquad \lim_{p\to \infty} Y_3^2(p) = -\infty, 
\end{align}
thus the curve
  $\Fix_{3+}^{2+}$ is a graph of a smooth function   $W_{\fix}^3(Y^2)$, $Y^2 \in (-\infty, + \infty).$ 

Let $p=p_3>0$ be the positive root of the equation   $p=2 \tanh p$, then $\ds Y_3^2(p_3) = 0, \  W_3^2(p_3)=\frac{\sinh^2 p_3-3}{6 \sqrt{\sinh p_3}}.$

Let $\sinh p_0 = 3$, then $p_0 = \ln (3 +\sqrt{10})$. It is known that there hold the inequalities  $\ln 2 <0.7$, $\ln 3< 1.1$. Now we estimate $p_0$:
\begin{align*}
\ln (6+\sqrt{10}-3)&<\ln 6+ \frac{\sqrt{10}-3}{6}< 0.7 +1.1+\frac{\sqrt{10}}{6}-\frac{1}{2} \\
									 &< \frac{13 \sqrt{10}}{30}+\frac{\sqrt{10}}{6}=   2 \tanh p_0,
\end{align*}
whence $p_0<2\tanh p_0$, then $p_3>p_0$, thus there holds the inequality $$\sinh p_3>3.$$ It follows that
$\ds W_3^2(p_3) > W_3^2(p_0) = \frac{1}{\sqrt{3}} > \frac{1}{\sqrt{\pi}}$,
i.e.,  $\ds W_{\fix}^3(0) > \frac{1}{\sqrt{\pi}}$. Whence it follows that 
\begin{align}
\exists \epsilon>0 \ \ \forall 0<Y^2<\epsilon \qquad -W_{\conj}^{22-}(-Y^2)<W_{\fix}^3(Y^2)<W_{\conj}^{21}(Y^2),
\end{align}
i.e., the statement of the lemma holds in a neighborhood of zero.

Now assume that there exists a point where this condition is false, i.e., 
there exists $Y_{\epsilon}^2 > 0$ such that $-W_{\conj}^{22-}(-Y_{\epsilon}^2)=W_{\conj}^{21}(Y_{\epsilon}^2)$. Then by continuity there is a point   $(Y^2,W^2) \in \Fix_{3+}^{2+} \cap (\mathcal{CI}_{x+}^{+} \cup \mathcal{CN}_{x--}^{+})$, which contradicts to Lemma~\ref{cutfix}. Consequently, our assumption is false, q.e.d. 
\end{proofn}

\begin{lem}\label{subset22}
There holds the inclusion $\Exp (\MAX_{2--}^{2+}) \subset \mathcal{N}_{x-}^+$.
\end{lem}
\begin{proofn}
 $\Exp (\MAX_{2--}^{2+})$ lies in the right half-plane of the plane   $(Y^2,W^2)$. 

Notice that for 
  $k\in (0,1), u_2 \in (\pi/2,\pi)$ there hold the conditions  
\begin{align}
J_{11}^{22}(k,u_2)&=-\frac{\partial Y_2^2(k,u_2)}{\partial u_2} = -\frac{\sqrt {\iota_4(k)} k^2 c_2 s_2}{2 \sqrt{d_2^3 \sqrt{1 - k^2}}} > 0, \label{J1122}\\
J_{21}^{22}(k,u_2)&=-\frac{\partial Y_2^2(k,u_2)}{\partial k} = \frac{k \Big((1 - k^2) \iota_6(k) + c_2^2 \big(\iota_4(k) + k^2 \iota_6(k)\big)\Big)}{2 \sqrt{d_2^3 \iota_4(k) \sqrt{(1 - k^2)^5}}}>0, \label{J2122} \\
\nabla^{22}(k,u_2) &= \frac{\partial W_2^2\big(k,u_2\big)/\partial k}{\partial Y_2^2 \big(k,u_2\big)/\partial k} - \frac{\partial W_2^2\big(k,u_2\big)/\partial u_2}{\partial Y_2^2 \big(k,u_2\big)/\partial u_2}  \nonumber \\
&=\frac{d_2^3 k^2 (E^2(k) - (1 - k^2) K^2(k))}{\iota_4^2(k) \sqrt{1 - k^2} \Big((1 - k^2) \iota_6(k) + c_2^2 \big(\iota_4(k) + k^2 \iota_6(k)\big)\Big)}>0, \label{cond22}
\end{align}
moreover $J_{11}^{22}(k,\pi/2) = 0$, $J_{21}^{22}(k,\pi/2)>0$, and $\nabla^{22}(k,\pi/2)>0$. Whence it follows from Lemma~\ref{below} and definition of $\mathcal{N}_{z-}^+$ that the set   $\Exp(\MAX_{--}^{1+})$ lies below the curve   $-W_{\conj}^{22-} (-Y^2), Y^2 \in (0, \infty)$, which defines the upper boundary of the set   $\mathcal{N}_{x-}^+$.

Now we should show that the set
 $\Exp(\MAX_{2--}^{2+})$ lies above the curve   $-W_{\conj}^{22+} (-Y^2), Y^2 \in (0, \infty)$. To this end consider in the preimage of the exponential mapping a symmetric set  $\MAX_{2++}^{2+}$ with a symmetric curve   $W_{\conj}^{22+} (Y^2)$, $Y^2 \in (-\infty, 0)$, with the inverted parameter   $\tilde{k} = 1 - k$, and $u_2 \in (0, \pi/2)$. In this case hypotheses of Lemma~\ref{below} hold as well, thus the set   $\MAX_{2++}^{2+}$ lies below the curve   $W_{\conj}^{22+} (Y^2), Y^2 \in (-\infty, 0)$. So the set    $\Exp(\MAX_{2--}^{2+})$ lies above the curve   $-W_{\conj}^{22+} (-Y^2), Y^2 \in (0, \infty)$. Consequently, $\Exp (\MAX_{2--}^{2+}) \subset \mathcal{N}_{x-}^+$.
\end{proofn}

\begin{proofn}[of Lemma~\ref{diffxm2}]
Apply Th.~\ref{Hadamard}. Conditions (1), (2) follow from definitions of the sets in the image and preimage of the exponential mapping. Condition (3) follows from results of~\cite{engel_conj}. By virtue of Lemma~\ref{subset22}, to prove condition (4) we consider arbitrary sequences $(k_n, u_{2}^n)$ in the image of  $\Exp$ tending to boundary of the set   $\MAX_{2--}^{2+}$, and show that such sequences, under the action of the exponential mapping, tend to boundary of  $\mathcal{N}_{z-}^+$: 
\begin{enumerate}
\item $k\to 0$, then $Y^2(k_n, u_{2}^n) \to 0$.
\item $k\to 1$, then $Y^2(k_n, u_{2}^n) \to \infty$.
\item $u_{2}^n \to \pi/2$, then by definition we have  $$\big(Y^2(k_n, u_{2}^n),W^2(k_n, u_{2}^n)\big) \to \CMAX_{2--}^{2+}.$$
\item $u_{2}^n \to \pi$, then by definition we have $$\big(Y^2(k_n, u_{2}^n),W^2(k_n, u_{2}^n)\big) \to \CMAX_{2-+}^{2+}.$$
\end{enumerate}
Whence by Th.~\ref{Hadamard} the mapping $\Exp (\MAX_{2--}^{2+}) \to \mathcal{N}_{x-}^+$ is a diffeomorphism.
\end{proofn}

\begin{proofn}[of Lemma~\ref{lem:smoothx2}]
With account of the symmetry
  $\varepsilon^4$, it suffices to consider the case    $\mathcal{CN}_{x+}^{+}$.

Lemma~\ref{lemcon22+} and Corollary~\ref{con22-} imply that  $\mathcal{CN}_{x+}^{+}$ form a decreasing continuous curve. Since the curves   $\mathcal{CN}_{x++}^{+}$ and $\mathcal{CN}_{x--}^{+}$ are continuous, in order to prove smoothness of union of these curves at the point   $\mathcal{CC}_{x+}^{+}$ it suffices to show that limits of the corresponding derivatives at this point coincide one with another. To this end evaluate the corresponding Taylor polynomials at the point  $k=0$:  
\begin{align*}
 \big(Y_2^2(k,0),W_2^2(k,0)\big) &= \bigg(-\frac{\sqrt{\pi} k^2}{4}+o(k^3), \frac{1}{\sqrt{\pi}} + \frac{3 k^2}{16 \sqrt{\pi}} +o(k^3)\bigg),\\
 \big(-Y_2^2(k,\pi/2),-W_2^2(k,\pi/2)\big) &= \bigg(\frac{\sqrt{\pi} k^2}{4}+o(k^3), \frac{1}{\sqrt{\pi}} - \frac{3 k^2}{16 \sqrt{\pi}} +o(k^3)\bigg).
\end{align*}
Now it follows that the curves
  $\mathcal{CN}_{x++}^{+}$ and $\mathcal{CN}_{x--}^{+}$ join one another smoothly at the point  $\mathcal{CC}_{x+}^{+}$. 
\end{proofn}

\section*{ACKNOWLEDGMENT}
This work was supported by the Russian Science Foundation 
under grant 17-11-01387 and performed in Ailamazyan Program Systems Institute 
of Russian Academy of Sciences.

\end{document}